\documentclass{article}%
\usepackage{latexsym}
\usepackage{amsmath}
\usepackage{graphicx}
\usepackage{amsfonts}
\usepackage{amssymb}%
\newtheorem{theorem}{Theorem}
\newtheorem{acknowledgement}[theorem]{Acknowledgement}

\newtheorem{corollary}[theorem]{Corollary}

\newtheorem{definition}[theorem]{Definition}

\newtheorem{lemma}[theorem]{Lemma}
\newtheorem{notation}[theorem]{Notation}
\newtheorem{problem}[theorem]{Problem}
\newtheorem{proposition}[theorem]{Proposition}
\newtheorem{remark}[theorem]{Remark}

\begin{document}

\author{Jamey Bass\\University of California\\Department of Mathematics,\\Santa Cruz, CA 95064
\and Andrey Todorov\thanks{Partially supported by The Institute of Mathematical
Sciences of The Chinese University of Hong Kong.}\\University of California\\Department of Mathematics,\\Santa Cruz, CA 95064\\Bulgarian Academy of Sciences\\Institute of Mathematics\\Sofia, Bulgaria}
\title{The Analogue of the Dedekind Eta Function for CY Manifolds I.}
\maketitle

\begin{abstract}
This is the first of series of articles in which we are going to study the
regularized determinants of the Laplacians of Calabi Yau metrics acting on
$(0,q)$ forms on the moduli space of CY manifolds with a fixed polarization.

It is well known that in the case of elliptic curves, the Kronecker limit
formula gives an explicit formula for the regularized determinant of the flat
metric with fixed volume on the elliptic curves. The following formula holds
in this case
\begin{equation}
\exp\left(  \frac{d}{ds}\zeta(s)\left\vert s=0\right.  \right)  =\det
\Delta_{(0,1)}(\tau)=\left(  \operatorname{Im}\tau\right)  ^{2}|\eta
(\tau)|^{4}, \label{10}%
\end{equation}
where $\eta(\tau)$ is the Dedekind eta function. It is well known fact that
$\eta(\tau)^{24}$ is a cusp automorphic form of weight $12$ related to the
discriminant of the elliptic curve.

Formula $\left(  \ref{10}\right)  $ implies that there exists a holomorphic
section of some power of the line bundle of the classes of cohomologies of
$(1,0)$ forms of the elliptic curves over the moduli space with an $L^{2}$
norm equal to $\det\Delta_{(0,1)}(\tau).$ This section is $\eta(\tau)^{24}.$
The purpose of this project is to find the relations between the regularized
determinants of CY metric on $(0,1)$ forms and algebraic discriminants of CY manifolds.

In this paper we will establish the local analogue of the formula $\left(
\ref{10}\right)  $ for CY manifolds.

\end{abstract}
\tableofcontents

\section{Introduction}

\subsection{General Comments}

The regularized determinants play an important role in different branches of
Physics and Mathematics. The explicit formulas for the regularized
determinants have many important applications both in Mathematics and Physics.
They are well known. This is the first of a series of papers in which we will
study the regularized determinant of the Laplacians of CY metrics on $(0,q)$
forms. We will first review the computation and the main results about the the
regularized determinant of the flat metric on the elliptic. We will describe
its relations to the discriminant locus and then we will describe our program
how to generalize the results in the case of elliptic curves to CY manifolds.

An easy computations show that the zeta function of the Laplace operator of
the flat metric on an elliptic curve with periods $(1,\tau)$ is given by the expression%

\[
E(s)=\left(  \frac{1}{2\pi}\right)  ^{2s}\underset{n,m\in\mathbb{Z}%
\text{\&}(n,m)\neq(0,0)}{\sum}\frac{1}{|n+m\tau|^{2s}}%
\]
where $\tau\in\mathbb{C},\operatorname{Im}\tau>0$ and ' means that the sum is
over all pair of integers $(m,n)\neq(0,0).$ The computation of the regularized
determinant in the case of the flat metric on an elliptic curve is based on
the Kronecker limit formula. See \cite{RS}. It states that $E(s)$ has a
meromorphic continuation in $\mathbb{C}$ with only one pole at $s=1$ and
\begin{equation}
\exp\left(  -\frac{d}{ds}E(s)|_{s=0}\right)  =\frac{1}{\left(
\operatorname{Im}\tau\right)  ^{2}\left\vert \eta\right\vert ^{4}} \label{KLF}%
\end{equation}
where $\eta$ is the Dedekind eta function. Thus formula $\left(
\ref{KLF}\right)  $ implies that
\begin{equation}
\det\Delta_{\tau,1}=\exp\left(  \frac{d}{ds}E(s)|_{s=0}\right)  =\left(
\operatorname{Im}\tau\right)  ^{2}\left\vert \eta\right\vert ^{4}.
\label{KLF0}%
\end{equation}
Let us denote by $G_{k}(\tau)$ the Eisenstein series of the lattice spanned by
$1$ and $\tau.$ Then the elliptic curve $E_{\tau}=\mathbb{C}/\left(
n+m\tau\right)  $ is given by
\[
y^{2}=4x^{3}-g_{2}x-g_{3},
\]
where $g_{2}=60G_{2}$ and $g_{3}=140G_{3}.$ The formula for the discriminant
$\Delta(\tau)$ of the elliptic curve is given by
\begin{equation}
\Delta(\tau)=\frac{g_{2}^{3}(\tau)-27g_{3}^{2}(\tau)}{1728g_{2}^{3}(\tau)}.
\label{8}%
\end{equation}
It is a well known formula that $\Delta(\tau)=\eta(\tau)^{24}.$ Thus the
formula for the regularized determinant of the Laplacian is related to the
discriminant of the elliptic curve $E_{\tau}.$

In this paper we will use the following method of deriving $\left(
\ref{KLF0}\right)  .$ See \cite{JT95} and \cite{JT96}. Compute the Hessian of
$\log\det\Delta_{\tau,1}$ and show that%
\begin{equation}
dd^{c}\log\det\Delta_{\tau,1}=-dd^{c}\log\operatorname{Im}\tau=dd^{c}%
\log\left\langle \omega_{\tau},\omega_{\tau}\right\rangle . \label{KLF1}%
\end{equation}
Thus $\left(  \ref{KLF1}\right)  $ implies that $\log\det\Delta_{\tau,1}$ is a
potential of the Poincare metric on the upper half plane $\mathfrak{h},$ which
is the Teichm\"{u}ller space of marked elliptic curves. Notice that%
\[
\operatorname{Im}\tau=-\frac{\sqrt{-1}}{2}%
{\displaystyle\int\limits_{E}}
\omega_{\tau}\wedge\overline{\omega_{\tau}}=\left\langle \omega_{\tau}%
,\omega_{\tau}\right\rangle ,
\]
where $\omega_{\tau}$ are holomorphic one forms on the elliptic curves
normalized as follows:%
\begin{equation}%
{\displaystyle\int\limits_{\gamma_{0}}}
\omega_{\tau}=1, \label{Nor}%
\end{equation}
where $\gamma_{0}$ is one of the generators of $H_{1}(E,\mathbb{Z}).$ This
means that: From the relation $\left(  \ref{KLF1}\right)  $ we derive that we
have
\begin{equation}
\det\Delta_{\tau,1}=\left\langle \omega_{\tau},\omega_{\tau}\right\rangle
\left\vert \eta(\tau)\right\vert ^{4}=\left(  \operatorname{Im}\tau\right)
^{2}\left\vert \eta\right\vert ^{4}, \label{KLF2}%
\end{equation}
where $\eta(\tau)$ is a holomorphic function.

We will outline how $\left(  \ref{KLF2}\right)  $ implies that $\eta(\tau)$ is
the Dedekind eta function. This can be done in several steps. The first step
is to prove that $\det\Delta_{\tau,1}$ is a bounded function on the moduli
space of the elliptic curves. The next step is to prove that $\left\langle
\omega_{\tau},\omega_{\tau}\right\rangle =\operatorname{Im}\tau$ has a
logarithmic growth near the the infinity of $\mathbb{PSL}_{2}(\mathbb{Z}%
)\backslash\mathfrak{h}.$ This will imply that $\left\vert \eta(\tau
)\right\vert $ must vanish at infinity. It is not difficult to see that
\begin{equation}
\mathbb{SL}_{2}(\mathbb{Z})/\left[  \mathbb{SL}_{2}(\mathbb{Z}),\mathbb{SL}%
_{2}(\mathbb{Z})\right]  \approxeq\mathbb{Z}/12\mathbb{Z}. \label{com}%
\end{equation}
The normalization $\left(  \ref{Nor}\right)  $ implies that $\left\langle
\omega_{\tau},\omega_{\tau}\right\rangle =\operatorname{Im}\tau$ is an
automorphic form of weight $-2.$ Thus $\left(  \ref{KLF2}\right)  $ and
$\left(  \ref{com}\right)  $ imply that $\eta(\tau)^{24}$ will be a cusp form
of weight 12 up to a constant. This will prove that $\eta(\tau)$ is the
Dedekind eta function. This fact intuitively represents the following
observation, when the metric degenerates then the discrete spectrum of the
Laplacian in the limit turns into a continuous spectrum which contains the
zero. Thus the regularized determinant must vanish at the points of the
compactified moduli space that corresponds to singular manifolds. $\det
\Delta_{\tau,1}$ is a bounded function on $\mathbb{PSL}_{2}(\mathbb{Z}%
)\backslash\mathfrak{h}.$

One of the observation in this paper is that the normalization of the period
of the holomorphic form $\omega_{\tau}$ $\left(  \ref{Nor}\right)  $ is
related to the choice of a maximal unipotent element of the mapping class
group as follows. It was proved in \cite{LTY} that maximal unipotent elements
in the mapping class group correspond to a monodromy operator $T$ acting on
the middle cohomology group of a generic fibre of a family of projective
algebraic varieties over the unit disk with only one singular fibre over the
origin of the disk. Grothendieck proved that we always have:
\[
(T^{N}-id)^{n+1}=0
\]
where $N$ is some positive integer and $n$ is the complex dimension of the
fibre of the family. Maximal unipotent elements corresponds to Jordan blocks
of dimension one. In case of CY manifolds it was proved in \cite{LTY} that if
$T$ has maximal index of unipotency then $T$ has a unique Jordan block of
dimension $n+1.$ Thus once we choose the unipotent element of the mapping
class group we can associate to it a unique up to a sign primitive class of
homology in the middle homology group which corresponds to the invariant
vanishing cycle. The cycle $\gamma_{0}$ that appears in $\left(
\ref{Nor}\right)  $ is the vanishing invariant cycle associated with the
maximal unipotent element.

It is a classical fact that the Dedekind eta function is related to the
algebraic discriminant of the equation that defines the elliptic curve. Thus
the arguments that we provided to prove that the absolute value of the
holomorphic function that appeared in $\left(  \ref{KLF2}\right)  $ is in fact
the analytic analogue of the discriminant. This follows from the fact that
$\det\Delta_{\tau,1}$ can be interpreted as the $L^{2}$ norm of a holomorphic
section of some power of the line bundle of the holomorphic one forms over the
moduli space of elliptic curves. In our next publications we will generalized
these arguments for CY manifolds.

This paper is the first one in the realization of the above described program.
In it we will compute $dd^{c}\log\det\Delta_{\tau,1}.$ We will show that
locally in the Teichm\"{u}ller space of the polarized CY manifolds we have%
\begin{equation}
dd^{c}\log\det\Delta_{\tau,1}=dd^{c}\log\left\langle \omega_{\tau}%
,\omega_{\tau}\right\rangle , \label{com1}%
\end{equation}
where $\omega_{\tau}$ is a family of holomorphic $n$ forms. Thus $\left(
\ref{com1}\right)  $ implies
\begin{equation}
\det\Delta_{\tau,1}=\left\langle \omega_{\tau},\omega_{\tau}\right\rangle
|\eta(\tau)|^{2}. \label{com2}%
\end{equation}
In fact by proving $\left(  \ref{com2}\right)  $ we generalized $\left(
\ref{KLF2}\right)  $ for higher dimensional CY manifolds.

In the case of elliptic curves formula $\left(  \ref{KLF1}\right)  $ implies
that $\log\det\Delta_{\tau,1}$ is the potential of the Poincare metric on the
upper half plane which is the Teichm\"{u}ller space of elliptic curves. In
\cite{To89} we proved that $dd^{c}\log\left\langle \omega_{\tau},\omega_{\tau
}\right\rangle $ gives the imaginary part of the Weil-Petersson metric on the
Teichm\"{u}ller space of polarized CY manifolds. Thus by establishing formula
$\left(  \ref{com1}\right)  ,$ we proved that the logarithm of the regularized
determinant $\det\Delta_{\tau,1}$ forms is the potential of the Weil Petersson metric.

\subsection{Outline of the Main Ideas}

The ideas that we used in this paper are similar to the ideas used in
\cite{DK} to established variational formulas for regularized determinants on
vector bundles on Riemann surfaces. We need first to fix the coordinates of
the local deformation space of the fix CY manifold M$_{0}$ with a fix
polarization class. In fact we use the coordinates $(\tau^{1},...,\tau^{N})$
introduced in \cite{To89}. These coordinates were later used in string theory.
See \cite{BCOV}. The coordinates $(\tau^{1},...,\tau^{N})$ depend on the
choice of a basis of harmonic forms $\phi_{1},...,\phi_{N}$ of $H^{1}($%
M$_{0},$T$^{1,0}).$ In these coordinates we have the following local
expression for the family $\overline{\partial_{\tau}}$ on M$_{0}:$%
\[
\overline{\partial_{i,\tau}}=\frac{\overline{\partial}}{\overline{\partial
z^{i}}}-%
{\displaystyle\sum\limits_{k=1}^{N}}
\tau^{k}\phi_{k,\overline{i}}^{j}\frac{\partial}{\partial z^{j}}+O(\tau^{2}),
\]
where $(z^{1},...,z^{n})$ are local coordinates on M$_{0}.$ Thus we get that
$\overline{\partial_{\tau}}$ depends holomorphically on $\tau$ and%
\begin{equation}
\frac{\partial}{\partial\tau^{i}}\overline{\partial_{\tau}}=-\phi_{i}%
\circ\partial_{0}+O(\tau). \label{2}%
\end{equation}
From the uniqueness of the solution of Calabi conjecture by Yau we know that
once we fix the polarization class then we obtain a family of Calabi Yau
metrics $g_{\tau}$ such that $\operatorname{Im}g_{\tau}=L.$ One of our results
is that $\operatorname{Im}g_{\tau}$ is a constant symplectic form. From here
and the expression for $\overline{\partial_{\tau}}^{\ast}$ we obtain that
$\overline{\partial_{\tau}}^{\ast}$ depends anti holomorphically on $\tau.$

After these preliminary results we are ready to compute the Hessian of
$\log\det\Delta_{\tau,q}.$ If we know the expressions for $\frac{\partial^{2}%
}{\partial\tau^{i}\overline{\partial\tau^{j}}}\log\Delta_{\tau,q}^{^{"}},$
where $\Delta_{\tau,q}^{^{"}}=\Delta_{\tau,q}|_{\operatorname{Im}%
\overline{\partial}^{\ast}}$ then we know how to compute $\frac{\partial^{2}%
}{\partial\tau^{i}\overline{\partial\tau^{j}}}\log\Delta_{\tau,q}.$ So we will
describe the ideas of the computations of $\frac{\partial^{2}}{\partial
\tau^{i}\overline{\partial\tau^{j}}}\log\Delta_{\tau,q}^{^{"}}.$ Since by
definition of the regularized determinants of the Laplacians are given by the
formula $\Delta_{\tau,q}^{^{"}}=\exp\left(  -\frac{d}{ds}\zeta_{\Delta
_{\tau,q}^{^{"}}}(s)\right)  |_{s=0}$ all our computations are based on the
following expression
\begin{equation}
\zeta_{\Delta_{\tau,q}^{^{"}}}(s)=\frac{1}{\Gamma(s)}%
{\displaystyle\int\limits_{0}^{\infty}}
Tr\left(  \frac{1}{4\pi t^{n}}\exp\left(  -t\Delta_{\tau,q}^{"}\right)
\right)  t^{s-1}dt. \label{3}%
\end{equation}
From $\left(  \ref{2}\right)  $ and $\left(  \ref{3}\right)  $ we obtain%
\[
\frac{\overline{\partial}}{\overline{\partial\tau^{i}}}\left(  \zeta_{q,\tau
}^{^{"}}(s)\right)  |_{\tau=0}=
\]
\begin{equation}
\frac{1}{\Gamma(s)}\int\limits_{0}^{\infty}Tr\left(  \exp(-t(\triangle
_{0,q}^{^{"}})\circ\triangle_{0,q}^{^{"}}\circ\partial_{0}^{-1}\circ
\mathcal{F}^{^{\prime}}(q+1,\overline{\phi_{i}})\circ\overline{\partial_{0}%
}\right)  t^{s}dt, \label{4}%
\end{equation}
where $\phi_{i}\in\mathbb{H}^{1}($M,$\Theta$) are the Kodaira-Spencer classes
viewed as bundle maps:%
\[
\phi_{i}:C^{\infty}\left(  \text{M,}\Omega_{\text{M}}^{1,0}\right)
\rightarrow C^{\infty}\left(  \text{M,}\Omega_{\text{M}}^{0,1}\right)
\]
and $\mathcal{F}\left(  q+1,\phi_{i}\right)  $ are the map induced by
\[
\phi_{i}\wedge id_{q-1}:C^{\infty}\left(  \text{M,}\Omega_{\text{M}}%
^{0,q}\right)  \rightarrow C^{\infty}\left(  \text{M,}\Omega_{\text{M}}%
^{1,0}\otimes\Omega_{\text{M}}^{0,q-1}\right)  .
\]
on $C^{\infty}\left(  \text{M,}\Omega_{\text{M}}^{0,q+1}\right)  $ and
restricted in $\operatorname{Im}\overline{\partial}.$ From $\left(
\ref{4}\right)  $ we obtain
\[
\frac{\partial^{2}}{\overline{\partial\tau^{j}}\partial\tau^{i}}\left(
\zeta_{\tau,q}^{"}(s)\right)  =
\]%
\[
\frac{1}{\Gamma(s)}\int\limits_{0}^{\infty}Tr\left(  \left(  \frac{\partial
}{\partial\tau^{j}}\exp(-t\left(  \triangle_{\tau,q}^{^{"}}\right)  \right)
\circ\left(  \triangle_{\tau,q}^{^{"}}\right)  \circ\partial_{\tau}^{-1}\circ
F^{^{\prime}}(q+1,\overline{\phi_{i}})\circ\overline{\partial_{0}}\right)
t^{s}dt+
\]%
\begin{equation}
\frac{1}{\Gamma(s)}\int\limits_{0}^{\infty}Tr\left(  \exp(-t\left(
\triangle_{\tau,q}^{^{"}}\right)  \circ\frac{\partial}{\partial\tau^{j}%
}\left(  \triangle_{\tau,q}^{^{"}}\right)  \circ\partial_{\tau}^{-1}\circ
F^{^{\prime}}(q+1,\overline{\phi_{i}})\circ\overline{\partial_{0}}\right)
t^{s}dt. \label{5}%
\end{equation}
Following the ideas in \cite{DK} to integrate by parts $\left(  \ref{5}%
\right)  $ we derive that:%
\begin{equation}
\frac{\partial^{2}}{\partial\tau^{i}\overline{\partial\tau^{j}}}\zeta_{\tau
,q}^{"}(s)=\frac{s(1-s)}{\Gamma(s)}%
{\displaystyle\int\limits_{0}^{\infty}}
Tr\left(  \exp\left(  -t\Delta_{\tau,q+1}^{^{\prime}}\right)  \mathcal{F}%
\left(  q+1,\phi_{i}\circ\overline{\phi_{j}}\right)  \right)  t^{s-1}dt.
\label{6}%
\end{equation}
By using the short term asymptotic expansion of $\exp\left(  -t\Delta
_{\tau,q+1}^{^{\prime}}\right)  $ we derive that
\begin{equation}
Tr\left(  \exp\left(  -t\Delta_{\tau,q+1}^{^{\prime}}\right)  \mathcal{F}%
\left(  q+1,\phi_{i}\circ\overline{\phi_{j}}\right)  \right)  =%
{\displaystyle\sum\limits_{k=-n}^{1}}
\frac{\alpha_{q,k}}{t^{k}}+\alpha_{q,0}+\psi_{q}(t). \label{7}%
\end{equation}
Combining $\left(  \ref{6}\right)  $ and $\left(  \ref{7}\right)  $ we get
that
\[
\frac{d}{ds}\left(  \frac{\partial^{2}}{\partial\tau^{i}\overline{\partial
\tau^{j}}}\zeta_{\tau,q}^{"}(s)\right)  |_{s=0}=\frac{\partial^{2}}%
{\partial\tau^{i}\overline{\partial\tau^{j}}}\log\det\Delta_{\tau,q}^{^{"}%
}=\alpha_{0,q}.
\]
In this paper we gave an explicit formula for $\alpha_{0,q}.$ By using the
expressions for $\alpha_{0,1}$ and $\alpha_{0,n-1}$ we show that $d\circ
d^{c}\left(  \log\det\Delta_{\tau,1}\right)  =-\operatorname{Im}WP,$ where
$\operatorname{Im}WP$ is the imaginary part of the Weil-Petersson metric. This
can be viewed as the generalization of the results in \cite{JT95}
and\ \cite{JT98}.

\subsection{Organization of the Paper}

This article is organized as follows.

In \textbf{Section 2} we introduce some basic notions about zeta functions of
Laplacians on Riemannian manifolds. We review the results from \cite{To89}.

In \textbf{Section 3 }We review the basic properties of the Weil-Petersson
metric on the moduli of CY manifolds. See also \cite{LS}.

In \textbf{Section 4 }we review the theory of moduli of CY manifolds following
\cite{LTYZ} and also metrics with logarithmic singularities on vector bundle
following Mumford's article \cite{Mu}. We prove that the $L^{2}$ metric on the
dualizing line bundle over the moduli space of CY manifolds is a good metric
in the sense of Mumford. This implies that the volumes of the moduli spaces of
CY manifolds are rational numbers.

In \textbf{Section 5} we review some facts about the Hilbert spaces of the
$(0,q)$ forms and their isospectral identifications which we used in the
paper. We also study traces the operators acting on the $L^{2}$ sections of
some vector bundle induced by some global $C^{\infty}$ section its
endormorphisms composed with the heat kernel. We study the short term
expansions of these operators and especially the constant term of the expansion.

In \textbf{Section 6 }we establish the variational formulas for the zeta
functions of the Laplacians and its regularized determinants. We also showed
that
\[
dd^{c}\log\det\Delta_{\tau,1}=-\operatorname{Im}W.P.
\]

In \textbf{Section 7} we prove that for K\"{a}hler Ricci flat manifold N with
$H^{0}($N,$\Omega_{\text{M}}^{n})=0$ we have
\[
dd^{c}\log\det\Delta_{\tau,n}=-\operatorname{Im}W.P.
\]

\begin{acknowledgement}
The first author wants to thank G. Moore, J. Jorgenson S. Donaldson, S.-T.
Yau, G. Zuckerman, D. Kazhdan, P. Deligne, S. Lang, B. Lian, J. Li, K. Liu, Y.
Eliashberg, Dai and R. Donagi for useful comments and support. I want to thank
Sinan Unver for his help. Special thanks to The Institute of Mathematical
Sciences at CUHK and Hangzhou for their hospitality during the preparation of
this article.
\end{acknowledgement}

\section{Preliminary Material}

\subsection{Basic Notions}

Let M be a n-dimensional K\"{a}hler manifold with a zero canonical class.
Suppose that $H^{k}($M$,\mathcal{O}_{\text{M}})=0$ for $1\leq k<n.$
\textit{Such manifolds are called Calabi-Yau manifolds. }A pair (M,$L$) will
be called a polarized CY manifold if M is a CY manifold and $L\in H^{2}%
($M,$\mathbb{Z}$)\footnote{Notice that $H^{1,1}($M,$\mathbb{R)=}H^{2}%
($M,$\mathbb{R)}$ since $H^{2}$(M,$\mathcal{O}_{\text{M}})=0$ for CY
manifolds.} is a fixed class such that it represents the imaginary part of a
K\"{a}hler metric on M.

Yau's celebrated theorem asserts the existence of a unique Ricci flat
K\"{a}hler metric g on M such that the cohomology class $[Im(g)]=L.$ (See
\cite{Yau}.) From now on we will consider polarized CY manifolds of odd
dimension. The polarization class $L$ determines the CY metric g uniquely. We
will denote by
\[
\bigtriangleup_{q}=\overline{\partial}^{\ast}\circ\overline{\partial
}+\overline{\partial}\circ\overline{\partial}^{\ast}%
\]
the associated Laplacians that act on smooth $(0,q)$ forms on M for $0\leq
q\leq n$. $\overline{\partial}^{\ast}$ is the adjoint operator of
$\overline{\partial}$ with respect to the CY metric g.

The regularized determinants are defined as follows: Let (M,g) be an
n-dimensional Riemannian manifold. Let
\[
\Delta_{q}=dd^{\ast}+d^{\ast}d
\]
be the Laplacian acting on the space of q forms on M. We recall that the
spectrum of the Laplacian $\Delta_{q}$ is positive and discrete. Thus the non
zero eigen values of $\Delta_{q}$ are
\[
0<\lambda_{1}\leq\lambda_{2}\leq...\leq\lambda_{n}\leq...
\]
We define the zeta function of $\Delta_{q}$ as follows:
\[
\zeta_{q}(s)=\sum_{i=1}^{\infty}\lambda_{i}^{-s}.
\]
It is known that $\zeta_{q}(s)$ is a well defined analytic function for
$\operatorname{Re}(s)\gg C,$ it has a meromorphic continuation in the complex
plane and $0$ is not a pole of $\zeta_{q}(s).$ Define
\[
\det(\Delta_{q})=\exp\left(  -\frac{d}{ds}\left(  \zeta_{q}(s)\right)
|_{s=0}\right)  .
\]
The determinant of $\ $these operators $\bigtriangleup_{q},$ defined through
zeta function regularization, will be denoted by det$\left(  \bigtriangleup
_{q}\right)  .$

The Hodge decomposition theorem asserts that
\[
\Gamma(\text{M},\Omega_{\text{M}}^{0,q})=\operatorname{Im}(\overline{\partial
})\oplus\operatorname{Im}(\overline{\partial}^{\ast})
\]
for $1\leq q\leq\dim_{\mathbb{C}}M-1.$ The restriction of $\ \bigtriangleup
_{q}$ on $\operatorname{Im}(\overline{\partial})$ will be denoted by
$\bigtriangleup_{q}^{^{\prime}}$ and $\ \bigtriangleup_{q}^{^{\prime}%
}=\overline{\partial}\circ\overline{\partial}^{\ast}$ and the restriction of
$\ \Delta_{q}$ on $\operatorname{Im}(\overline{\partial}^{\ast})$ will be
denoted by $\bigtriangleup_{q}^{"}$ and $=\overline{\partial}^{\ast}%
\circ\overline{\partial}.$ Hence we have
\[
Tr(\exp(-t\bigtriangleup_{q})=Tr(\exp(-t\bigtriangleup_{q}^{^{\prime}%
})+Tr(\exp(-t\bigtriangleup_{q}^{"}).
\]
This implies that
\[
\zeta_{q}(s)=\sum_{k=1}^{\infty}\lambda_{k}^{-s}=\zeta_{q}^{^{\prime}%
}(s)+\zeta_{q}^{"}(s),
\]
where $\lambda_{k}>0$ are the positive eigen values of $\bigtriangleup_{q}$
and $\zeta_{q}^{^{\prime}}(s)$ \& $\zeta_{q}^{"}(s)$ are the zeta functions of
$\bigtriangleup_{q}^{^{\prime}}$ and $\bigtriangleup_{q}^{"}.$ From here and
the definition of \ the regularized determinant we obtain that
\[
\log\det(\bigtriangleup_{q})=\log\det(\bigtriangleup_{q}^{^{\prime}})+\log
\det(\bigtriangleup_{q}^{"}).
\]
It is a well known fact that the action of $\bigtriangleup_{q}^{^{"}}$ on
$\operatorname{Im}\overline{\partial}^{\ast}$ is isospectral to the action of
$\bigtriangleup_{q+1}^{^{\prime}}$ on $\operatorname{Im}\overline{\partial},$
which means that the spectrum of $\bigtriangleup_{q}^{^{"}}$ is equal to the
spectrum of $\bigtriangleup_{q+1}^{^{\prime}}.$ So we have the equality
\[
\det(\bigtriangleup_{q}^{"})=\det(\bigtriangleup_{q+1}^{^{\prime}}).
\]
\ 

\begin{notation}
Let f be a map from a set A to a set B and let g be a map from the set B to
the set C, then the compositions of those two maps we will denote by f$\circ$g.
\end{notation}

\subsection{Basic Notions about Complex Structures}

Let M be an even dimensional C$^{\infty}$ manifold. We will say that M has an
almost complex structure if there exists a section $I\in C^{\infty}%
($M$,Hom(T^{\ast},T^{\ast})$ such that $I^{2}=-id.$ $T$ is the tangent bundle
\ and $T^{\ast}$ is the cotangent bundle on M. This definition is equivalent
to the following one: Let M be an even dimensional C$^{\infty}$ manifold.
Suppose that there exists a global splitting of the complexified cotangent
bundle
\[
T_{\text{M}}^{\ast}\otimes\mathbf{C}=\Omega_{\text{M}}^{1,0}\oplus
\Omega_{\text{M}}^{0,1},
\]
where $\Omega_{\text{M}}^{0,1}=\overline{\Omega_{\text{M}}^{1,0}}.$ Then we
will say that M has an almost complex structure. We will say that an almost
complex structure is an integrable one, if for each point $x\in$M there exists
an open set $U\subset$M such that we can find local coordinates $z^{1}%
,..,z^{n},$ such that $dz^{1},..,dz^{n}$ are linearly independent in each
point $m\in U$ and they generate $\Omega_{\text{M}}^{1,0}|_{U}.$

\begin{definition}
\label{belt}Let M be a complex manifold. Let $\phi\in\Gamma($M$,Hom(\Omega
_{\text{M}}^{1,0},\Omega_{\text{M}}^{0,1})),$ then we will call $\phi$ a
Beltrami differential.
\end{definition}

Since
\[
\Gamma(\text{M},Hom(\Omega_{\text{M}}^{1,0},\Omega_{\text{M}}^{0,1}%
))\backsimeq\Gamma(\text{M},\Omega_{\text{M}}^{0,1}\otimes T_{\text{M}}%
^{1,0}),
\]
we deduce that locally $\phi$ can be written as follows:
\[
\phi|_{U}=\sum\phi_{\overline{\alpha}}^{\beta}\overline{dz}^{\alpha}%
\otimes\frac{\partial}{\partial z^{\beta}}.
\]
From now on we will denote by $A_{\phi}$ the following linear operator:
\[
A_{\phi}=\left(
\begin{array}
[c]{cc}%
id & \phi(\tau)\\
\overline{\phi(\tau)} & id
\end{array}
\right)  .
\]
We will consider only those Beltrami differentials $\phi$ such that
det($A_{\phi})\neq0.$ The Beltrami differential \ $\phi$\ defines an
integrable complex structure on M if and only if the following equation
holds:
\begin{equation}
\overline{\partial}\phi=\frac{1}{2}\left[  \phi,\phi\right]  , \label{0}%
\end{equation}
where
\begin{equation}
\left[  \phi,\phi\right]  |_{U}:=\sum_{\nu=1}^{n}\sum_{1\leqq\alpha<\beta\leqq
n}\left(  \sum_{\mu=1}^{n}\left(  \phi_{\overline{\alpha}}^{\mu}\left(
\partial_{\mu}\phi_{\overline{\beta}}^{\nu}\right)  -\phi_{\overline{\beta}%
}^{\mu}\left(  \partial_{\nu}\phi_{\overline{\alpha}}^{\nu}\right)  \right)
\right)  \overline{dz}^{\alpha}\wedge\overline{dz}^{\beta}\otimes
\frac{\partial}{dz^{\nu}} \label{eq:def}%
\end{equation}
(See \cite{KM}.)

\subsection{Kuranishi Space and Flat Local Coordinates}

Kuranishi proved the following Theorem:

\begin{theorem}
\label{Kur}Let $\left\{  \phi_{i}\right\}  $ be a basis of harmonic $(0,1)$
forms of $\mathbb{H}^{1}($M$,T^{1,0})$ on a Hermitian manifold M. Let $G$ be
the Green operator and let $\phi(\tau^{1},..,\tau^{N})$ be defined as
follows:\
\begin{equation}
\phi(\tau^{1},...,\tau^{N})=\sum_{i=1}^{N}\phi_{i}\tau^{i}+\frac{1}%
{2}\overline{\partial}^{\ast}G[\phi(\tau^{1},...,\tau^{N}),\phi(\tau
^{1},...,\tau^{N})]. \label{eq:bel}%
\end{equation}
\textit{There exists }$\varepsilon>0$\textit{\ such that if }$\tau=(\tau
^{1},...,\tau^{N})$ \textit{satisfies }$|\tau_{i}|<\varepsilon$\textit{ then
}$\phi(\tau^{1},...,\tau^{N})$\textit{\ is a\ global }$C^{\infty}%
$\textit{\ section of the bundle }$\Omega_{\text{M}}^{(0,1)}\otimes T^{1,0}%
$.(See \cite{KM}.)
\end{theorem}

Based on Theorem \ref{Kur}, we proved in \cite{To89} the following Theorem:

\begin{theorem}
\label{tod1}Let M be a CY manifold and let $\left\{  \phi_{i}\right\}  $ be a
basis of harmonic $(0,1)$ forms with coefficients in $T^{1,0},$ i.e.
\[
\left\{  \phi_{i}\right\}  \in\mathbb{H}^{1}(\text{M},T^{1,0}),
\]
then the equation $\left(  \ref{0}\right)  $ has a solution in the form:
\[
\phi(\tau^{1},...,\tau^{N})=\sum_{i=1}^{N}\phi_{i}\tau^{i}+\sum_{|I_{N}%
|\geqq2}\phi_{I_{N}}\tau^{I_{N}}=
\]%
\[
\sum_{i=1}^{N}\phi_{i}\tau^{i}+\frac{1}{2}\overline{\partial}^{\ast}%
G[\phi(\tau^{1},...,\tau^{N}),\phi(\tau^{1},...,\tau^{N})]
\]
and $\overline{\partial}^{\ast}\phi(\tau^{1},...,\tau^{N})=0,$ $\phi_{I_{N}%
}\lrcorner\omega_{\text{M}}=\partial\psi_{I_{N}}$ \textit{where }
$I_{N}=(i_{1},...,i_{N})$\ \ \textit{is a multi-index}, \
\[
\phi_{I_{N}}\in C^{\infty}(\text{M},\Omega_{\text{M}}^{0,1}\otimes
T^{1,0}),\tau^{I_{N}}=(\tau^{i})^{i_{1}}...(\tau^{N})^{i_{N}}%
\]
\textit{and if for some} $\varepsilon>0$ $|\tau^{i}|<\varepsilon$ then
\textit{\ } $\phi(\tau)\in C^{\infty}($M$,\Omega_{\text{M}}^{0,1}\otimes
T^{1,0})$ \textit{where} $i=1,...,N.$ $\left(  See\text{ }\cite{Ti}\text{
}and\text{ }\cite{To89}\right)  .$
\end{theorem}

It is a standard fact from Kodaira-Spencer-Kuranishi deformation theory that
for each%
\[
\tau=(\tau^{1},...,\tau^{N})\in\mathcal{K}%
\]
as in Theorem \ref{tod1} the Beltrami differential $\phi(\tau^{1},...,\tau
^{N})$ defines a new integrable complex structure on M. This means that the
points of $\mathcal{K},$ where
\[
\mathcal{K}:\{\tau=(\tau^{1},...,\tau^{N})||\tau^{i}|<\varepsilon\}
\]
defines a family of operators $\overline{\partial}_{\tau}$ on the $C^{\infty}$
family
\[
\mathcal{K}\times\text{M}\rightarrow\text{M}%
\]
and $\overline{\partial}_{\tau}$ are integrable in the sense of
Newlander-Nirenberg. Moreover it was proved by Kodaira, Spencer and Kuranishi
that we get a complex analytic family of CY manifolds $\pi
:\mathcal{X\rightarrow K},$ where as $C^{\infty}$ manifold
$\mathcal{X\backsimeq K}\times$M$.$ The family
\begin{equation}
\pi:\mathcal{X\rightarrow K} \label{kur}%
\end{equation}
is called the Kuranishi family. The operators $\overline{\partial}_{\tau}$ are
defined as follows:

\begin{definition}
\label{tod3}Let $\{\mathcal{U}_{i}\}$ be an open covering of M, with local
coordinate system in $\mathcal{U}_{i}$ given by $\{z_{i}^{k}\}$ with
$k=1,...,n=$dim$_{\mathbb{C}}$M. Assume that $\phi(\tau^{1},...,\tau
^{N})|_{\mathcal{U}_{i}}$ is given by:
\[
\phi(\tau^{1},...,\tau^{N})=\sum_{j,k=1}^{n}(\phi(\tau^{1},...,\tau
^{N}))_{\overline{j}}^{k}\text{ }d\overline{z}^{j}\otimes\frac{\partial
}{\partial z^{k}}.
\]
\textit{Then we define \ }\
\begin{equation}
\left(  \overline{\partial}\right)  _{\tau,\overline{j}}=\frac{\overline
{\partial}}{\overline{\partial z^{j}}}-\sum_{k=1}^{n}(\phi(\tau^{1}%
,...,\tau^{N}))_{\overline{j}}^{k}\frac{\partial}{\partial z^{k}}.
\label{eq:dibar}%
\end{equation}

\end{definition}

\begin{definition}
\label{flat}The coordinates $\tau=(\tau^{1},...,\tau^{N})$ defined in Theorem
\ref{tod1}, will be fixed from now on and will be called the flat coordinate
system in $\mathcal{K}$.
\end{definition}

\subsection{Family of Holomorphic Forms}

In \cite{To89} the following Theorem is proved:

\begin{theorem}
\label{forms}There exists a family of holomorphic forms $\omega_{\tau}$ of the
Kuranishi family $\left(  \ref{kur}\right)  $ such that%
\[
\left\langle \lbrack\omega_{\tau}],[\omega_{\tau}]\right\rangle =
\]%
\[
1-%
{\displaystyle\sum\limits_{i,j}}
\left\langle \omega_{0}\lrcorner\phi_{i},\omega_{0}\lrcorner\phi
_{j}\right\rangle \tau^{i}\overline{\tau^{j}}+%
{\displaystyle\sum\limits_{i,j}}
\left\langle \omega_{0}\lrcorner\left(  \phi_{i}\wedge\phi_{k}\right)
,\omega_{0}\lrcorner\left(  \phi_{j}\wedge\phi_{l}\right)  \right\rangle
\tau^{i}\overline{\tau^{j}}\tau^{k}\overline{\tau^{l}}+O(\tau^{5})=
\]%
\[
1-%
{\displaystyle\sum\limits_{i,j}}
\tau^{i}\overline{\tau^{j}}+%
{\displaystyle\sum\limits_{i,j}}
\left\langle \omega_{0}\lrcorner\left(  \phi_{i}\wedge\phi_{k}\right)
,\omega_{0}\lrcorner\left(  \phi_{j}\wedge\phi_{l}\right)  \right\rangle
\tau^{i}\overline{\tau^{j}}\tau^{k}\overline{\tau^{l}}+O(\tau^{5})\text{ and}%
\]%
\begin{equation}
\left\langle \lbrack\omega_{\tau}],[\omega_{\tau}]\right\rangle \leq
\left\langle \lbrack\omega_{0}],[\omega_{0}]\right\rangle . \label{form}%
\end{equation}

\end{theorem}

\section{Weil-Petersson Metric}

\subsection{Basic Properties}

It is a well known fact from Kodaira-Spencer-Kuranishi theory that the tangent
space $T_{\tau,\mathcal{K}\text{ }}$at a point $\tau\in\mathcal{K}$ can be
identified with the space of harmonic (0,1) forms with values in the
holomorphic vector fields $\mathbb{H}^{1}($M$_{\tau},T$). We will view each
element $\phi\in\mathbb{H}^{1}($M$_{\tau},T$) as a point wise linear map from
$\Omega_{\text{M}_{\tau}}^{(1,0)}$ to $\Omega_{\text{M}_{\tau}}^{(0,1)}.$
Given $\phi_{1}$ and $\phi_{2}\in\mathbb{H}^{1}($M$_{\tau},T$)$,$ the trace of
the map
\[
\phi_{1}\circ\overline{\phi_{2}}:\Omega_{\text{M}_{\tau}}^{(0,1)}%
\rightarrow\Omega_{\text{M}_{\tau}}^{(0,1)}%
\]
at the point $m\in$M$_{\tau}$ with respect to the metric g is
simply:\textit{\ \ }\
\begin{equation}
Tr(\phi_{1}\circ\overline{\phi_{2}})=\sum_{k,l,m=1}^{n}(\phi_{1}%
)_{\overline{l}}^{k}(\overline{\phi_{2})_{\overline{k}}^{m}}g^{\overline{l}%
,k}g_{k,\overline{m}} \label{wp}%
\end{equation}

\begin{definition}
\label{WP}We will define the Weil-Petersson metric on $\mathcal{K}$ via the
scalar product:\textit{\ \ }\
\begin{equation}
\left\langle \phi_{1},\phi_{2}\right\rangle =\int\limits_{\text{M}}Tr(\phi
_{1}\circ\overline{\phi_{2}})vol(g). \label{wp1}%
\end{equation}

\end{definition}

We proved in \cite{To89} that the coordinates
\[
\tau=(\tau^{1},...,\tau^{N})
\]
as defined in Definition \ref{flat} are flat in the sense that the
Weil-Petersson metric is K\"{a}hler and in these coordinates we have that the
components $g_{i,\overline{j}}$ of the Weil Petersson metric are given by the
following formulas:%

\[
g_{i,\overline{j}}=\delta_{i,\overline{j}}+R_{i,\overline{j},l,\overline{k}%
}\tau^{l}\overline{\tau^{k}}+O(\tau^{3}).
\]

Very detailed treatment of the Weil-Petersson geometry of the moduli space of
polarized CY manifolds can be found in \cite{Lu} and \cite{LS}. In those two
papers important results are obtained.

\subsection{Infinitesimal Deformation of the Imaginary Part of the WP Metric}

\begin{theorem}
\label{const}Near each point $\tau_{0}$ of the Kuranishi space $\mathcal{K},$
the imaginary part $\operatorname{Im}(g)$ of the CY metric $g$ has the
following expansion in the coordinates $\tau:=(\tau^{1},...,\tau^{N})$:
\[
\operatorname{Im}(g)(\tau,\overline{\tau})=\operatorname{Im}(g)(\tau
_{0})+O((\tau-\tau_{0})^{2}).
\]

\end{theorem}

\textbf{Proof: }Without loss of generality we can assume that $\tau_{0}=0.$ In
\cite{To89} we proved that the forms\
\begin{equation}
\theta_{\tau}^{k}=dz^{k}+\sum_{l=1}\left(  \phi(\tau^{1},...,\tau
^{N})_{\overline{l}}^{k}\right)  d\overline{z^{l}}\text{ } \label{eq:m}%
\end{equation}
for $k=1,...,n$ form a basis of $(1,0)$ forms relative to the complex
structure defined by $\tau\in\mathcal{K}$ in $\mathcal{U\subset}$M. Let\ \
\begin{equation}
\operatorname{Im}(g_{\tau})=\sqrt{-1}\left(  \sum_{1\leq k\leq l\leq
n}g_{k,\overline{l}}(\tau,\overline{\tau})\theta_{\tau}^{k}\wedge
\overline{\theta_{\tau}^{l}}\right)  \label{eq:m0}%
\end{equation}
and\
\begin{equation}
g_{k,\overline{l}}(\tau,\overline{\tau})=g_{k,\overline{l}}(0)+\sum_{i=1}%
^{N}\left(  \left(  g_{k,\overline{l}}(1)\right)  _{i}\tau^{i}+\left(
g_{k,\overline{l}}^{^{\prime}}(1)\right)  _{\overline{i}}\overline{\tau^{i}%
}\right)  +O(2). \label{WP1a}%
\end{equation}
We get the following expression for $\operatorname{Im}(g_{\tau})$ in terms of
$dz^{i}$ and $\overline{dz^{j}}$, by substituting the expressions for
$\theta_{\tau}^{k}$ from $\left(  \ref{eq:m}\right)  $ and the expressions for
$g_{k,\overline{l}}(\tau,\overline{\tau})$ from formula $\left(
\ref{WP1a}\right)  $ in the formula $\left(  \ref{eq:m0}\right)  $:%

\[
\operatorname{Im}(g_{\tau})=\sqrt{-1}\left(  \sum_{1\leq k\leq l\leq
n}g_{k,\overline{l}}(\tau,\overline{\tau})\theta_{\tau}^{k}\wedge
\overline{\theta_{\tau}^{l}}\right)  =\sqrt{-1}\left(  \sum_{1\leq k\leq l\leq
n}g_{k,\overline{l}}(0)dz^{k}\wedge\overline{dz^{l}}\right)  +
\]

\[
+\sqrt{-1}\left(  \sum_{i=1}^{N}\tau^{i}\left(  \sum_{1\leq k\leq l\leq
n}\left(  \left(  g_{k,\overline{l}}(1)\right)  _{i}dz^{k}\wedge
\overline{dz^{l}}+\sum_{m=1}^{n}(g_{k,\overline{m}}\overline{\phi
_{i,\overline{l}}^{m}}-g_{l,\overline{m}}\overline{\phi_{i,\overline{k}}^{m}%
})dz^{k}\wedge dz^{l}\right)  \right)  \right)
\]

\[
+\frac{1}{\sqrt{-1}}\sum_{i=1}^{N}\overline{\tau^{i}}\overline{\left(
\sum_{1\leq k\leq l\leq n}\left(  \left(  g_{k,\overline{l}}(1)\right)
_{i}dz^{k}\wedge\overline{dz^{l}}+(\sum_{m=1}^{n}(g_{k,\overline{m}}%
\overline{\phi_{i,\overline{l}}^{m}}-g_{l,\overline{m}}\overline
{\phi_{i,\overline{k}}^{m}})dz^{k}\wedge dz^{l}\right)  \right)  }.
\]
On page 332 of \cite{To89} the following results is proved:

\begin{lemma}
\label{sym}Let $\phi\in\mathbb{H}^{1}($M$,T$) be a harmonic form with respect
to the CY metric g. Let
\[
\phi|_{U}=\sum_{k,l=1}^{n}\phi_{\overline{k}}^{l}\text{ }\overline{dz}%
^{k}\otimes\frac{\partial}{\partial z^{l}},
\]
\textit{then }
\[
\phi_{\overline{k},\overline{l}}=\sum_{j=1}^{n}g_{j,\overline{k}}\text{ }%
\phi_{\overline{l}}^{j}=\sum_{j=1}^{n}g_{j,\overline{l}}\text{ }%
\phi_{\overline{k}}^{j}=\phi_{\overline{l},\overline{k}}.
\]

\end{lemma}

From Lemma \ref{sym} we conclude that
\begin{equation}
\sum_{m=1}^{n}(g_{k,\overline{m}}\overline{\phi_{i,\overline{l}}^{m}%
}-g_{l,\overline{m}}\overline{\phi_{i,\overline{k}}^{m}})=0. \label{Eq}%
\end{equation}
From $\left(  \ref{Eq}\right)  $ we get the following expression for
$\operatorname{Im}(g_{\tau})$:
\[
\operatorname{Im}(g_{\tau})=\sqrt{-1}\left(  \sum_{1\leq k\leq l\leq
n}g_{k,\overline{l}}(0)dz^{k}\wedge\overline{dz^{l}}\right)  +
\]%
\[
\sqrt{-1}\left(  \sum_{i=1}^{N}\tau^{i}\left(  \sum_{1\leq k\leq l\leq
n}\left(  g_{k,\overline{l}}(1)\right)  _{i}dz^{k}\wedge\overline{dz^{l}%
}\right)  \right)  +
\]%
\begin{equation}
+\sqrt{-1}\left(  \sum_{i=1}^{N}\overline{\tau^{i}}\overline{\sum_{1\leq k\leq
l\leq n}\left(  g_{k,\overline{l}}(1)\right)  _{i}dz^{k}\wedge\overline
{dz^{l}}}\right)  +O(2) \label{eq:m3}%
\end{equation}
Let us define the (1,1) forms\textit{\ }$\psi_{i}$ as follows:\textit{\ \ }\
\begin{equation}
\psi_{i}=\sqrt{-1}\left(  \sum_{1\leq k\leq l\leq n}\left(  g_{k,\overline{l}%
}(1)\right)  _{i}dz^{k}\wedge\overline{dz^{l}}\right)  \label{eq:m4}%
\end{equation}
We derive the following formula, by substituting in the expression $\left(
\ref{eq:m3}\right)  $ the expression given by $\left(  \ref{eq:m4}\right)
$:\textit{\ \ }\
\begin{equation}
\operatorname{Im}(g_{\tau})=\operatorname{Im}(g_{0})+\sum_{i=1}^{N}\tau
^{i}\psi_{i}+\sum_{i=1}^{N}\overline{\tau^{i}\psi_{i}}+O(\tau^{2}) \label{m7}%
\end{equation}
From the fact that the class of the cohomology of the imaginary part of the CY
metric is fixed, i.e.\textit{\ } $[\operatorname{Im}(g_{\tau}%
)]=[\operatorname{Im}(g_{0})]=L,$ and $\left(  \ref{m7}\right)  $ we deduce
that\textit{\ }each\textit{\ }$\psi_{i}$\textit{\ }is an exact form\textit{,
i.e. \ }\
\begin{equation}
\psi_{i}=\sqrt{-1}\partial\overline{\partial}f_{i}, \label{m5}%
\end{equation}
where \textit{\ }$f_{i}$ are globally defined functions on\textit{\ }%
M\textit{. }Our Theorem will follow if we prove that \textit{\ }$\psi_{i}=0.$

\begin{lemma}
\label{const1}$\psi_{i}=0.$
\end{lemma}

\textbf{Proof: }In \textit{\cite{To89}\ }we proved that\textit{\ }
\begin{equation}
\det(g_{\tau})=\wedge^{n}\operatorname{Im}(g_{\tau})=\det(g_{0})+O(2)
\label{WP0}%
\end{equation}
in the flat coordinates $(\tau^{1},...,\tau^{N}).$ We deduce from the
expressions $\left(  \ref{m7}\right)  $ and $\left(  \ref{m5}\right)  $, by
direct computations that:%

\[
\det(g_{\tau})=\det(g_{0})+\sqrt{-1}\sum_{i=1}^{N}\tau^{i}\left(  \sum
_{k,l}g^{\overline{l},k}\partial_{k}\overline{\partial_{l}}(f_{i})\right)  +
\]%
\begin{equation}
\frac{1}{\sqrt{-1}}\sum_{i=1}^{N}\overline{\tau^{i}}\overline{\left(
\sum_{k,l}g^{\overline{l},k}\partial_{k}\overline{\partial_{l}}(f_{i})\right)
}+O(2). \label{WP1}%
\end{equation}
Combining $\left(  \ref{WP0}\right)  $ and $\left(  \ref{WP1}\right)  $ we
obtain that for each i we have\textit{: }
\[
\sum_{k,l}g^{\overline{l},k}\partial_{k}\overline{\partial_{l}}(f_{i}%
)=\triangle(f_{i})=0,
\]
where $\triangle$ is the Laplacian of the metric g. From the maximum
principle, we deduce that all\textit{\ }$f_{i}$ are constants. Formula
$\left(  \ref{m5}\right)  $ implies that $\psi_{i}=0.$ Lemma \ref{const1} is
proved. $\blacksquare$

Theorem \ref{const} follows directly from Lemma \ref{const1}. Theorem
\ref{const} is proved\textit{. }$\blacksquare$

\begin{corollary}
\label{L2a}The imaginary part $\operatorname{Im}$g$_{\tau}$ of the CY metric
is a constant symplectic form on the moduli space $\mathfrak{M}_{L}$(M).
\end{corollary}

\begin{corollary}
\label{L2} \ The following formulas are true:
\begin{equation}
\frac{\partial}{\partial\tau_{i}}\left(  \overline{\partial_{\tau}}\right)
^{\ast}=0\text{ and }\frac{\overline{\partial}}{\overline{\partial\tau_{i}}%
}\left(  \partial_{\tau}\right)  ^{\ast}=0. \label{v5}%
\end{equation}

\end{corollary}

\textbf{Proof: }We know from K\"{a}hler geometry that $(\overline
{\partial_{\tau}})^{\ast}=[\Lambda_{\tau},\partial_{\tau}],$ where
$\Lambda_{\tau}$ is the contraction with \textit{(1,1) }vector field$:$%
\begin{equation}
\Lambda_{\tau}=\frac{\sqrt{-1}}{2}\sum_{k,l=1}^{n}g_{\tau}^{\overline{k}%
,l}(\theta_{\tau}^{l})^{\ast}\wedge(\overline{\theta_{\tau}^{k}})^{\ast}
\label{v6}%
\end{equation}
on M$_{\tau}$ and $(\theta_{\tau}^{l})^{\ast}$ is (1,0) vector field on
M$_{\tau}$ dual to the (1,0) form
\[
\theta_{\tau}^{i}=dz^{i}+\sum_{j=1}^{N}\tau^{j}(\sum_{k=1}^{n}(\phi
_{j})_{\overline{k}}^{i}\text{ }\overline{dz}^{k})).
\]
Corollary \ref{L2a}\textit{\ }implies that $\frac{\partial}{\partial\tau_{i}%
}(\Lambda_{\tau})=0.$ On the other hand,\textit{\ }$\partial_{\tau}%
$\textit{\ }depends antiholomorphically on $\tau$, i.e. it depends
on\textit{\ }$\overline{\tau}=(\overline{\tau_{1}},...,\overline{\tau_{N}}).$
So we deduce that\textit{:}%

\[
\frac{\partial}{\partial\tau_{i}}((\overline{\partial_{\tau}})^{\ast})=\left(
\left[  \frac{\partial}{\partial\tau_{i}}(\Lambda_{\tau}),\partial_{\tau
}\right]  +\left[  \Lambda_{\tau},\frac{\partial}{\partial\tau_{i}}%
(\partial_{\tau})\right]  \right)  =0.
\]
Exactly in the same way we prove that $\frac{\overline{\partial}}%
{\overline{\partial\tau_{i}}}\left(  \partial_{\tau}\right)  ^{\ast}=0.$
Corollary \ref{L2} is proved. $\blacksquare$

\section{Moduli of CY\ Manifolds}

\subsection{Basic Construction}

\begin{definition}
\label{Teich}We will define the Teichm\"{u}ller space $\mathcal{T}$(M) of a CY
manifold M as follows: $\mathcal{T}($M$):=\mathcal{I}($M$)/Diff_{0}($M$),$
\textit{where}\
\[
\mathcal{I}(\text{M}):=\left\{  \text{all integrable complex structures on
M}\right\}
\]
\textit{and } Diff$_{0}$(M) \textit{is the group of diffeomorphisms isotopic
to identity. The action of the group Diff(M}$_{0})$ \textit{is defined as
follows; Let }$\phi\in$Diff$_{0}$(M) \textit{then }$\phi$ \textit{acts on
integrable complex structures on M by pull back, i.e. if }
\[
I\in C^{\infty}(M,Hom(T(\text{M}),T(\text{M})),
\]
\textit{then we define } $\phi(I_{\tau})=\phi^{\ast}(I_{\tau}).$
\end{definition}

\begin{definition}
We will call a pair $($M$;\gamma_{1},...,\gamma_{b_{n}})$ a marked CY manifold
where M is a CY manifold and $\{\gamma_{1},...,\gamma_{b_{n}}\}$ is a basis of
$H_{n}$(M,$\mathbb{Z}$)/Tor.
\end{definition}

\begin{remark}
\label{mark}Let $\mathcal{K}$ be the Kuranishi space. It is easy to see that
if we choose a basis of $H_{n}$(M,$\mathbb{Z}$)/Tor in one of the fibres of
the Kuranishi family $\pi:\mathcal{M\rightarrow K}$ then all the fibres will
be marked, since as a $C^{\infty}$ manifold $\mathcal{X}_{\mathcal{K}%
}\approxeq$M$\times\mathcal{K}$.
\end{remark}

In \cite{LTYZ} the following Theorem was proved:

\begin{theorem}
\label{teich}There exists a family of marked polarized CY manifolds
\begin{equation}
\mathcal{Z}_{L}\mathcal{\rightarrow}\widetilde{T}(\text{M}), \label{fam2}%
\end{equation}
which possesses the following properties: \textbf{a)} It is effectively
parametrized, \textbf{b) }The base has dimension $h^{n-1,1},$ \textbf{c) }For
any marked CY manifold M of fixed topological type for which the polarization
class $L$ defines an imbedding into a projective space $\mathbb{CP}^{N},$
there exists an isomorphism of it (as a marked CY manifold) with a fibre
M$_{s}$ of the family $\mathcal{Z}_{L}.$
\end{theorem}

\begin{corollary}
\label{teich1}Let $\mathcal{Y\rightarrow}$X be any family of marked polarized
CY manifolds, then there exists a unique holomorphic map $\phi:$%
X$\rightarrow\widetilde{T}($M$)$ up to a biholomorphic map $\psi$ of M which
induces the identity map on $H_{n}($M$,\mathbb{Z}).$
\end{corollary}

From now on we will denote by $\mathcal{T}$(M) the irreducible component of
the Teichm\"{u}ller space that contains our fixed CY manifold M.

\begin{definition}
We will define the mapping class group $\Gamma_{1}($M$)$ of any compact
C$^{\infty}$ manifold M as follows: $\Gamma_{1}$(M$)=Diff_{+}($M$)/Diff_{0}%
($M$),$ where $Diff_{+}($M$)$ is the group of diffeomorphisms of M preserving
the orientation of M and $Diff_{0}($M$)$ is the group of diffeomorphisms
isotopic to identity.
\end{definition}

\begin{definition}
Let $L\in H^{2}($M$,\mathbb{Z})$ be the imaginary part of a K\"{a}hler metric.
Let $\Gamma_{2}:=\{\phi\in\Gamma_{1}($M$)|\phi(L)=L\}.$
\end{definition}

It is a well know fact that the moduli space of polarized algebraic manifolds
$\mathcal{M}_{L}($M$)=\mathcal{T}($M$)/\Gamma_{2}.$ In \cite{LTYZ} the
following fact was established:

\begin{theorem}
\label{Vie}There exists a subgroup of finite index $\Gamma_{L}$ of
$\ \Gamma_{2}$ such that $\Gamma_{L}$ acts freely on $\mathcal{T}$(M) and
$\Gamma\backslash\mathcal{T}$(M)$=\mathfrak{M}_{L}$(M) is a non-singular
quasi-projective variety. Over $\mathfrak{M}_{L}$(M) there exists a family of
polarized CY manifolds $\mathcal{M}\rightarrow\mathfrak{M}_{L}($M$).$
\end{theorem}

\begin{remark}
\label{Vie1}Theorem \ref{Vie} implies that we constructed a family of
non-singular CY manifolds $\pi:\mathcal{X\rightarrow}\mathfrak{M}_{L}($M$)$
over a quasi-projective non-singular variety $\mathfrak{M}_{L}($M$)$. Moreover
it is easy to see that $\mathcal{X\subset}\mathbb{CP}^{N}\times\mathfrak{M}%
_{L}($M$).$ \textit{So} $\mathcal{X}$ \ \textit{is also quasi-projective. From
now on we will work only with this family.}
\end{remark}

\section{Hilbert Spaces and Trace Class Operators}

\subsection{Preliminary Material}

\begin{definition}
\label{Hilb1}We will denote by $L_{0,q}^{2}(\operatorname{Im}(\overline
{\partial}^{\ast}))$ the Hilbert subspace in $L^{2}($M,$\Omega_{\text{M}%
}^{(0,q)})$ which is the $L^{2}$ completion of \ $\overline{\partial^{\ast}}$
exact forms in $C^{\infty}($M,$\Omega_{\text{M}}^{(0,q)})$ for $q\geq0.$ In
the same manner we will denote by $L_{1,q-1}^{2}(\operatorname{Im}(\partial))$
the Hilbert subspace in $L^{2}($M,$\Omega_{\text{M}}^{(1,q-1)})$ which is the
$L^{2}$ competition of the $\partial$ exact $(1,q-1)$ forms in $C^{\infty}%
($M,$\Omega_{\text{M}}^{(1,q-1)})$ for $q>0$ . All the completions are with
respect to the scalar product on the bundles $\Omega_{\text{M}}^{p,q}$ defined
by the CY metric g.
\end{definition}

Let $\phi(\tau^{1},..,\tau^{N}$ ) be a solution of the equation $\left(
\ref{0}\right)  $:
\[
\overline{\partial}\phi(\tau^{1},...,\tau^{N})=\frac{1}{2}[\phi(\tau
^{1},...,\tau^{N}),\phi(\tau^{1},...,\tau^{N})]
\]
which is guaranteed by Theorem \ref{tod1}. From Definition \ref{belt} of the
Beltrami differential, we know that the Beltrami differential $\phi(\tau
^{1},..,\tau^{N}$ ) defines a linear fibrewise map
\[
\phi(\tau^{1},...,\tau^{N}):\Omega_{\text{M}}^{(1,0)}\rightarrow
\Omega_{\text{M}}^{(0,1)}.
\]
So\textit{\ \ }\
\begin{equation}
\phi(\tau^{1},...,\tau^{N})\in C^{\infty}(\text{M},Hom(\Omega_{\text{M}%
}^{(1,0)},\Omega_{\text{M}}^{(0,1)}). \label{tr1}%
\end{equation}

\begin{definition}
\label{Hilb}We define the following maps between vector bundles
\[
\phi\wedge id:\Omega_{\text{M}}^{(1,q-1)}\rightarrow\Omega_{\text{M}}^{(0,q)}%
\]
as
\[
\phi(dz^{i}\wedge\alpha)=\phi(dz^{i})\wedge\alpha
\]
\textit{for each }$1\leq q\leq n.$ Clearly each fibre wise linear map
$\phi\wedge id_{q-1}$ \textit{defines a natural linear operator }
\[
F(q,\phi):L^{2}(\text{M},\Omega_{\text{M}}^{(1,q-1)})\rightarrow
L^{2}(\text{M},\Omega_{\text{M}}^{(0,q)})
\]
between the Hilbert spaces. The restriction of \ the linear operator
F(q,$\phi$) on the subspace $\operatorname{Im}(\partial)\subset L^{2}%
($M$,\Omega_{\text{M}}^{(1,q-1)}))$ to $\operatorname{Im}(\overline{\partial
})\subset L^{2}($M$,\Omega_{\text{M}}^{(0,q)})$ will be denoted by $F^{\prime
}(q,\phi)$. The restriction of \ the linear operator F(q,$\phi$) on the
subspace $\operatorname{Im}(\partial^{\ast})\subset L^{2}($M$,\Omega
_{\text{M}}^{(1,q-1)}))$ to $\operatorname{Im}\left(  \overline{\partial
}^{\ast}\right)  \subset L^{2}($M$,\Omega_{\text{M}}^{(0,q)})$ will be denoted
by $F^{"}(q,\phi)$. Let $\phi$ and $\psi$ be two Kodaira Spencer classes and
let
\[
\phi\circ\overline{\psi}:L^{2}(\text{M},\Omega_{\text{M}}^{(0,1)})\rightarrow
L^{2}(\text{M},\Omega_{\text{M}}^{(0,1)})
\]
be fibrewise linear map given by
\begin{equation}
\phi\circ\overline{\psi}|_{U}:=\sum_{\alpha,\beta=1}^{n}\left(  \phi
\circ\overline{\psi}\right)  _{\overline{\beta}}^{\overline{\alpha}}%
\overline{dz^{\beta}}\otimes\frac{\partial}{\overline{\partial z^{\alpha}}}.
\label{hilb}%
\end{equation}
We define the fibrewise bundle maps
\begin{equation}
\left(  \phi\circ\overline{\psi}\right)  \wedge id_{q-1}:\Omega_{\text{M}%
}^{0,q}\rightarrow\Omega_{\text{M}}^{0,q}, \label{v3q}%
\end{equation}
as follows:
\begin{equation}
\left(  \left(  \phi\circ\overline{\psi}\right)  \wedge id_{q-1}\right)
(\omega):=\left(  \sum_{\alpha,\beta=1}^{n}\left(  \phi\circ\overline{\psi
}\right)  _{\overline{\beta}}^{\overline{\alpha}}\overline{dz^{\beta}}%
\otimes\frac{\overline{\partial}}{\overline{\partial z^{\alpha}}}\right)
\lrcorner\omega, \label{LO}%
\end{equation}
where $\lrcorner$ means contraction of tensors and $\omega$ is some global
form of type $(0,q)$ on M. We will define for the linear operators
\begin{equation}
\mathcal{F}^{^{\prime}}(q,\phi\circ\overline{\psi}):L_{0,q}^{2}%
(\operatorname{Im}(\overline{\partial}))\rightarrow L_{0,q}^{2}%
(\operatorname{Im}(\overline{\partial})) \label{v3a}%
\end{equation}
and
\begin{equation}
\mathcal{F}^{"}(q,\phi\circ\overline{\psi}):L_{0,q}^{2}(\operatorname{Im}%
(\overline{\partial})^{\ast})\rightarrow L_{0,q}^{2}(\operatorname{Im}%
(\overline{\partial})^{\ast}) \label{v4a}%
\end{equation}
as the restriction of the operators $\left(  \left(  \phi\circ\overline{\psi
}\right)  \wedge id_{q-1}\right)  $ on $L_{0,q}^{2}(\operatorname{Im}%
(\overline{\partial}))$ and $L_{0,q}^{2}(\operatorname{Im}(\overline{\partial
})^{\ast})$ respectively.
\end{definition}

\begin{remark}
\label{?}It is a standard fact that we can choose globally $\overline
{\partial}$ closed forms $\omega_{1},...,\omega_{N}$ of type $(0,q)$ such that
at each point $z\in$M they span the fibre $\Omega_{\text{M,}z}^{0,q}.$ We can
deduce directly from the definitions of the operators $F^{\prime}%
(q,\phi),\mathcal{F}^{^{\prime}}(q,\phi\circ\overline{\psi})$ and
$F^{^{\prime}}(q,\overline{\psi}\circ\phi)$ and the existence of the forms
$\omega_{1},...,\omega_{N}$ that the operators $F^{\prime}(q,\phi
),\mathcal{F}^{^{\prime}}(q,\phi\circ\overline{\psi})$ and $F^{^{\prime}%
}(q,\overline{\psi}\circ\phi)$ pointwise will be represented by matrices of
dimensions $\binom{n}{q},\binom{n}{q}$ and $n\times\binom{n}{q-1}.$
\end{remark}

\subsection{Trace Class Operators (See \cite{BGV})}

Let $H$ be a Hilbert space with a orthonormal basis $e_{i}.$ An operator $A$
is a \textbf{Hilbert-Schmidt} operator if%
\[
\left\Vert A\right\Vert _{HS}^{2}=%
{\displaystyle\sum\limits_{i}}
\left\Vert Ae_{i}\right\Vert ^{2}=%
{\displaystyle\sum\limits_{ij}}
\left\vert \left\langle Ae_{i},e_{j}\right\rangle \right\vert ^{2}<\infty
\]
is finite. The number $\left\Vert A\right\Vert _{HS}^{2}$ is called the
Hilbert-Schmidt norm of A. If A is a Hilbert-Schmid so is its adjoint
$A^{\ast}$ and $\left\Vert A\right\Vert _{HS}^{2}=\left\Vert A^{\ast
}\right\Vert _{HS}^{2}.$ If $U$ is a bounded operator on $H$ and $A$ is an
Hilbert-Schmidt, then $U\circ A$ and $A\circ U$ are Hilbert-Schmidt operators
and $\left\Vert U\circ A\right\Vert _{HS}\leq\left\Vert A\circ U\right\Vert
_{HS}.$

In this paper we will consider the Hilbert spaces of the square integrable
sections of the bundles $\Omega_{\text{M}}^{0,q}\otimes\left\vert
\Lambda_{\text{M}}\right\vert ^{1/2}$ on M, where $\left\vert \Lambda
_{\text{M}}\right\vert $ is the trivial density bundles generated by the
volume form of the CY metric.

An operator $K$ with square-integrable kernel%
\[
k(w,z)\in\Gamma_{L^{2}}\left(  \text{M}\times\text{M,}\left(  \Omega
_{\text{M}}^{0,q}\otimes\left\vert \Lambda_{\text{M}}\right\vert
^{1/2}\boxtimes\Omega_{\text{M}}^{0,q}\otimes\left\vert \Lambda_{\text{M}%
}\right\vert ^{1/2}\right)  \right)
\]
is Hilbert-Schmidt, and%
\begin{equation}
\left\Vert K\right\Vert _{HS}^{2}=%
{\displaystyle\int\limits_{(w,z)\in\text{M}\times\text{M}}}
Tr\left(  k\left(  w,z\right)  ^{\ast}k\left(  w,z\right)  \right)  .
\label{hs}%
\end{equation}
Formula $\left(  \ref{hs}\right)  $ follows from the definition of the
Hilbert-Schmidt norm
\[
\left\Vert K\right\Vert _{HS}^{2}=%
{\displaystyle\sum\limits_{i,j}}
\left\vert \left\langle Ke_{i},e_{j}\right\rangle \right\vert ^{2}.
\]
An operator $K$ is said to be \textbf{trace class }if it has the form $A\circ
B,$ where $A$ and $B$ are Hilbert-Schmidt. For such operators, the sum
\[
TrK=%
{\displaystyle\sum\limits_{i}}
\left\langle Ke_{i},e_{i}\right\rangle
\]
is absolutely summable and $TrK$ is independent of the choice of the
orthonormal basis in $H$ and is called the trace of $K.$

\subsection{Adiabatic Limits, Heat Kernels and Traces}

In this subsection we study the traces of operators which are compositions of
the heat kernel with operators induced by endomorphisms of some vector bundle.
We will use some of the results from \cite{DK} and will adopt them to our situation.

Let $h$ be a metric on a vector bundle $E$ over M. Let $\Delta_{h}$ be the
Laplacian on $E$. It is a well known fact that the operator $\exp
(-t\triangle_{h})$ can be represented by an integral kernel:%
\[
k_{t}(w,z,\tau)=%
{\displaystyle\sum\limits_{j}}
\exp\left(  -t\lambda_{j}\right)  \phi_{j}(w)\otimes\phi_{j}(z),
\]
where $\lambda_{j}$ and $\phi_{j}$ are the eigen values and the eigen sections
of the Laplace operator $\Delta_{h}$ on some vector bundle $E$ on M$.$
$k_{t}(w,z,\tau)$ is an operator of trace class. We know that the following
formula holds for the short term asymptotic expansion of $Tr\left(
k_{t}(w,z,\tau)\right)  $%
\[
Tr\left(  k_{t}(w,z,\tau)\right)  =\frac{\alpha_{-n}}{t^{n}}+...+\frac
{\alpha_{-1}}{t}+\alpha_{0}+O(t).
\]
Let $E$ be a holomorphic vector bundle over M, let $\phi\in C^{\infty}%
($M,Hom($E,E$)). It is easy to see that the operator $\exp(-\Delta_{h}%
)\circ\phi$ is of trace class and its trace has an asymptotic expansion
\begin{equation}
Tr\left(  k_{t}(w,z,\tau)\circ\phi\right)  =\frac{\beta_{-n}(\phi)}{t^{n}%
}+...+\frac{\beta_{-1}(\phi)}{t}+\beta_{0}(\phi)+O(t) \label{As}%
\end{equation}
according to \cite{BGV}. We will study the following problem in this section:

\begin{problem}
Find an explicit expression for $\beta_{0}(\phi).$
\end{problem}

\begin{definition}
\label{BGV1}We define the function $k_{\tau}^{d}(w,z,t)$ in a neighborhood of
the diagonal $\Delta$ in M$\times$M as follows: Let $\rho_{\tau}$ be the
injectivity radius on M$_{\tau}.$ Let $d_{\tau}(w,z)$ be the distance between
the points $w$ and $z$ on M$_{\tau}$ with respect to CY metric g$_{\tau}.$ We
suppose that $|\tau|<\varepsilon.$ Let $\delta$ be such that $\delta
>\rho_{\tau}.$ Then we define the function $k_{\tau}^{d}(w,z,t)$ as a
C$^{\infty}$ function on M$\times$M using partition of unity by using the
functions
\begin{equation}
k_{t}^{d}(w,z,\tau)=\left\{
\begin{array}
[c]{ll}%
\frac{1}{(4\pi t)^{\frac{n}{2}}}\exp\left(  -\frac{d_{\tau}^{2}(w,z)}%
{4t}\right)  & if\text{ }d_{\tau}(w,z)<\rho_{\tau}\\
0 & if\text{ }d_{\tau}(w,z)>\delta.
\end{array}
\right.  \label{bgv1}%
\end{equation}
defined on the opened balls around countable points $(w_{k},z_{k})$ on
M$\times$M with injectivity radius $\rho_{\tau}.$
\end{definition}

Let $E$ be a holomorphic vector bundle on M with a Hermitian metric $h$ on it
and let $\mathcal{P}_{\tau}(w,z)$ be the parallel transport of the bundle $E$
along the minimal geodesic joining the point $w$ and $z$ with respect to
natural connection on $E$ induced by the metric $h$ on $E.$ \ It was proved in
\cite{BGV} on page 87 that we can represent the operator $\exp(-t\triangle
_{h})$ by an integral kernel $k_{t}(w,z,\tau)$, where%

\begin{equation}
k_{t}(w,z,\tau)=k_{t}^{d}(w,z,\tau)\left(  \mathcal{P}_{\tau}%
(w,z)+O(t)\right)  \label{bgv}%
\end{equation}
and $\Delta_{h}:=\overline{\partial}_{h}^{\ast}\circ\overline{\partial}.$

\begin{definition}
\label{ker}We will define the kernel $k_{t}^{\#}(w,z,\tau)$ as the matrix
operator defined by
\begin{equation}
k_{t}^{\#}(w,z,\tau)=k_{t}^{d}(w,z,\tau)\mathcal{P}_{\tau}(w,z), \label{bgv2}%
\end{equation}
where
\begin{equation}
k_{t}(w,z,\tau)=k_{t}^{\#}(w,z,\tau)+\varepsilon_{t}(w,z,\tau). \label{bgv3}%
\end{equation}
Let us define
\begin{equation}
\Upsilon_{t}(\phi,\tau,z):=\int\limits_{\text{M}}Tr\left(  \left(  k_{t/2}%
^{d}(w,z,\tau)\right)  ^{\ast}\circ\left(  k_{t/2}^{\#}(w,z,\tau)\circ
\phi\right)  \right)  vol(g)_{w}. \label{me2}%
\end{equation}

\end{definition}

\begin{proposition}
\label{Trace0}We have%
\begin{equation}
\underset{t\rightarrow0}{\lim}\int\limits_{\text{M}}Tr\left(  \varepsilon
_{t}(w,z,0)\circ\phi\right)  vol(g)_{w}=0. \label{350}%
\end{equation}

\end{proposition}

\textbf{Proof: }The definition \ref{ker} of $k_{\tau}^{\#}(w,z,t)$ and the
arguments from \cite{DK} on page 260 imply that $\varepsilon_{t}(w,z,0)$ is
bounded and tends to zero away from the diagonal, as $t$ tends to zero. From
here we deduce that
\[
\underset{t\rightarrow0}{\lim}\int\limits_{\text{M}}Tr\left(  \varepsilon
_{t}(w,z,0)\circ\phi\right)  vol(g)_{w}=0
\]
uniformly in $z$. Proposition \ref{Trace0} is proved. $\blacksquare$

\begin{lemma}
\label{tr}Let $E$ be a holomorphic vector bundle over M, let $\phi\in
C^{\infty}($M,Hom($E,E$))$,$ then $\underset{t\rightarrow0}{\lim}\Upsilon
_{t}(\phi,\tau,z)$ exists and\textit{\ } \ \
\begin{equation}
\underset{t\rightarrow0}{\lim}\Upsilon_{t}(\phi,\tau,z)=Tr(\phi|_{E_{z}}).
\label{tr30}%
\end{equation}

\end{lemma}

\textbf{Proof:} We have:%

\[
\underset{t\rightarrow0}{\lim}\Upsilon_{t}(\phi,0,z)=\underset{t\rightarrow
0}{\lim}\int\limits_{\text{M}}Tr\left(  \left(  k_{t/2}^{d}(w,z,\tau)\right)
^{\ast}\circ\left(  k_{t/2}^{\#}(w,z,\tau)\circ\phi\right)  \right)
vol(g)_{w}=
\]

\begin{equation}
\underset{t\rightarrow0}{\lim}\int\limits_{\text{M}}Tr\left(  \frac{1}{(4\pi
t)^{\frac{n}{2}}}\exp\left(  -\frac{d_{0}^{2}(w,z)}{4t}\right)  \circ
\mathcal{P}_{0}(w,z)\circ\phi\right)  vol(g)_{w}. \label{expl}%
\end{equation}
Using the fact that%
\begin{equation}
\underset{t\rightarrow0}{\lim}\left(  \frac{1}{(4\pi t)^{\frac{n}{2}}}%
\exp\left(  -\frac{d_{\tau}^{2}(w,z)}{4t}\right)  \right)  =\delta(z-w),
\label{delta}%
\end{equation}%
\begin{equation}
\underset{w\rightarrow z}{\lim}\mathcal{P}_{0}(w,z)=id \label{delta1}%
\end{equation}
and the explicit formula $\left(  \ref{expl}\right)  $ for $\Upsilon_{t}%
(\phi,\tau,z)$ we obtain that\textit{\ } \ \
\begin{equation}
\underset{t\rightarrow0}{\lim}\Upsilon_{t}(\phi,\tau,z)=\int\limits_{\text{M}%
}Tr\left(  \delta(z-w)\circ\phi\right)  vol(g)_{\omega}=Tr(\phi(z)|_{E_{z}}).
\label{tr36}%
\end{equation}
Lemma \ref{tr} is proved. $\blacksquare$

\begin{theorem}
\label{Trace00}Let $\phi\in C^{\infty}\left(  \text{M},\left(  \Omega
_{\text{M}}^{0,q}\right)  ^{\ast}\otimes\Omega_{\text{M}}^{0,q}\right)  $ then
the operator $\exp(-t\Delta_{h})\circ\phi$ for $t>0$ is of trace class and its
trace is given by the formula;%
\begin{equation}
Tr\left(  \exp(-t\Delta_{h})\circ\phi\right)  =%
{\displaystyle\int\limits_{\text{M}}}
\Upsilon_{t}(\phi,\tau,z)vol(g)_{z}+\Phi(t), \label{me2a}%
\end{equation}
where the short term asymptotic of $\Phi(t)$ is given by
\begin{equation}
\Phi(t)=%
{\displaystyle\sum\limits_{k=1}^{N_{0}>0}}
\frac{a_{-k}}{t^{k}}+O(t). \label{me2b}%
\end{equation}

\end{theorem}

\textbf{Proof: }The proof of Theorem \ref{Trace00} is based on the facts that
\begin{equation}
\exp(-t\Delta_{h})\circ\phi=\exp\left(  -\frac{t}{2}\Delta_{h}\right)
\circ\exp\left(  -\frac{t}{2}\Delta_{h}\circ\phi\right)  , \label{me3}%
\end{equation}
and the operators $\exp\left(  -\frac{t}{2}\Delta_{h}\right)  $ and
$\exp\left(  -\frac{t}{2}\Delta_{h}\circ\phi\right)  $ can be represented by
$C^{\infty}$ kernels $k_{1}(z,w,t)$ and $k_{\phi}(z,w,t).$

As we pointed out the operators defined by the kernels $k_{1}(z,w,t)$ and
$k_{\phi}(z,w,t)$ are Hilbert-Schmidt operators. Thus since the operator
$\exp(-t\Delta_{h})\circ\phi$ is a product of two Hilbert-Schmidt operators it
is of trace class. On the other hand the definition of the trace of the
operator $\exp(-t\Delta_{h})\circ\phi$ implies that%
\[
Tr\left(  \exp(-t\Delta_{h})\circ\phi\right)  =\left\langle \exp\left(
-\frac{t}{2}\Delta_{h}\right)  ^{\ast},\exp\left(  -\frac{t}{2}\Delta
_{h}\right)  \circ\phi\right\rangle =
\]%
\begin{equation}%
{\displaystyle\int\limits_{(z,w)\in\text{M}\times\text{M}}}
Tr\left(  \left(  k_{1}(z,w,t)\right)  ^{\ast}\circ k_{\phi}(z,w,t)\right)
vol(g_{z,w}). \label{me4}%
\end{equation}
From the definitions of the function $\Upsilon_{t}(\phi,\tau,z)$ and the
operator $\varepsilon_{t}(w,z,\tau)$ we deduce that%
\[
Tr\left(  \exp(-t\Delta_{h})\circ\phi\right)  =\left\langle \left(
k_{t/2}(w,z,\tau)\right)  ,k_{t/2}(w,z,\tau)\circ\phi\right\rangle =
\]%
\[%
{\displaystyle\int\limits_{\text{M}}}
\left(
{\displaystyle\int\limits_{\text{M}}}
\left(  k_{t/2}(w,z,\tau)\right)  ^{\ast}\circ\left(  k_{t/2}(w,z,\tau
)\circ\phi\right)  vol(g)_{w}\right)  vol(g)_{z}=
\]%
\[%
{\displaystyle\int\limits_{\text{M}}}
\Upsilon_{t}(\phi,\tau,z)vol(g)_{z}+
\]%
\[%
{\displaystyle\int\limits_{\text{M}}}
\left(
{\displaystyle\int\limits_{\text{M}}}
Tr\left(  \left(  \varepsilon_{t/2}(w,z,\tau)\right)  ^{\ast}\circ
k_{t/2}^{\#}(w,z,\tau)\circ\phi\right)  vol(g)_{\omega}\right)  vol(g)_{z}+
\]%
\begin{equation}%
{\displaystyle\int\limits_{\text{M}}}
\left(
{\displaystyle\int\limits_{\text{M}}}
Tr\left(  k_{t/2}^{\#}(w,z,\tau)\circ\varepsilon_{t/2}(w,z,\tau)\circ
\phi\right)  vol(g)_{\omega}\right)  vol(g)_{z}. \label{me4a}%
\end{equation}

\begin{lemma}
\label{Trace1}Let
\[
\Phi_{1}(t):=%
{\displaystyle\int\limits_{\text{M}}}
\left(
{\displaystyle\int\limits_{\text{M}}}
Tr\left(  \left(  \varepsilon_{t/2}(w,z,\tau)\right)  ^{\ast}\circ
k_{t/2}^{\#}(w,z,\tau)\circ\phi\right)  vol(g)_{\omega}\right)  vol(g)_{z}%
\]
and
\begin{equation}
\Phi_{2}(t):=%
{\displaystyle\int\limits_{\text{M}}}
\left(
{\displaystyle\int\limits_{\text{M}}}
Tr\left(  k_{t/2}^{\#}(w,z,\tau)\circ\varepsilon_{t/2}(w,z,\tau)\circ
\phi\right)  vol(g)_{\omega}\right)  vol(g)_{z} \label{TR1}%
\end{equation}
then we have%
\begin{equation}
\Phi_{1}(t)=%
{\displaystyle\sum\limits_{k=1}^{N_{0}>0}}
\frac{b_{-k}}{t^{k}}+O(t)\text{ and }\Phi_{2}(t)=%
{\displaystyle\sum\limits_{k=1}^{N_{0}>0}}
\frac{c_{-k}}{t^{k}}+O(t). \label{TR2}%
\end{equation}

\end{lemma}

\textbf{Proof: }Let
\begin{equation}
k_{t/2}^{\#}(w,z,\tau)=%
{\displaystyle\sum\limits_{k=1}^{N_{0}>0}}
\frac{\mathcal{B}_{-k}(w,z)}{t^{k}}+\mathcal{B}_{0}(w,z)+%
{\displaystyle\sum\limits_{k=1}}
\mathcal{B}_{k}(w,z)t^{k} \label{TR3}%
\end{equation}
be the short term asymptotic expansion of the operator $k_{t/2}^{\#}%
(w,z,\tau).$ We know that
\begin{equation}
\underset{t\rightarrow0}{\lim}\varepsilon_{t}(w,z,\tau)=0 \label{TR4}%
\end{equation}
away from the diagonal $\Delta\subset$M$\times$M. Combining $\left(
\ref{TR3}\right)  $ and $\left(  \ref{TR4}\right)  $ with the definitions of
operators $k_{t/2}^{\#}(w,z,\tau)\circ\varepsilon_{t/2}(w,z,\tau)\circ\phi$
and $\left(  \varepsilon_{t/2}(w,z,\tau)\right)  ^{\ast}\circ k_{t/2}%
^{\#}(w,z,\tau)\circ\phi$ we obtain that%
\[
k_{t/2}^{\#}(w,z,\tau)\circ\varepsilon_{t/2}(w,z,\tau)\circ\phi=
\]%
\begin{equation}%
{\displaystyle\sum\limits_{k=1}^{N_{0}>0}}
\frac{\mathcal{B}_{-k}(w,z)\circ\varepsilon_{t/2}(w,z,\tau)\circ\phi}{t^{k}%
}+\mathcal{B}_{0}(w,z)\circ\varepsilon_{t/2}(w,z,\tau)\circ\phi+O(t)
\label{TR4a}%
\end{equation}
and
\[
\left(  \varepsilon_{t/2}(w,z,\tau)\right)  ^{\ast}\circ k_{t/2}^{\#}%
(w,z,\tau)\circ\phi=
\]%
\begin{equation}%
{\displaystyle\sum\limits_{k=1}^{N_{0}>0}}
\frac{\left(  \varepsilon_{t/2}(w,z,\tau)\right)  ^{\ast}\circ\mathcal{B}%
_{-k}(w,z)}{t^{k}}+\left(  \varepsilon_{t/2}(w,z,\tau)\right)  ^{\ast}%
\circ\mathcal{B}_{0}(w,z)+O(t). \label{TR4b}%
\end{equation}
Combining $\left(  \ref{TR4a}\right)  $, $\left(  \ref{TR4b}\right)  $ with
$\left(  \ref{TR4}\right)  $ we get that%
\[
\underset{t\rightarrow0}{\lim}\mathcal{B}_{0}(w,z)\circ\varepsilon
_{t/2}(w,z,\tau)\circ\phi=\underset{t\rightarrow0}{\lim}\left(  \varepsilon
_{t/2}(w,z,\tau)\right)  ^{\ast}\circ\mathcal{B}_{0}(w,z)=0
\]
away from the diagonal. From here we we obtain that
\[
\underset{t\rightarrow0}{\lim}%
{\displaystyle\int\limits_{\text{M}}}
Tr\left(  \mathcal{B}_{0}(w,z)\circ\varepsilon_{t/2}(w,z,\tau)\circ
\phi\right)  vol(g)=0
\]
and%
\[
\underset{t\rightarrow0}{\lim}%
{\displaystyle\int\limits_{\text{M}}}
Tr\left(  \left(  \varepsilon_{t/2}(w,z,\tau)\right)  ^{\ast}\circ
\mathcal{B}_{0}(w,z)\right)  vol(g)=0.
\]
Lemma \ref{Trace1} is proved. $\blacksquare$

Theorem \ref{Trace00} follows directly from Lemma \ref{Trace1} and $\left(
\ref{me4a}\right)  $. $\blacksquare$

\begin{theorem}
\label{Trace01}We have the following expression for $\beta_{0}(\phi)$ from
$\left(  \ref{As}\right)  $:%
\begin{equation}
\beta_{0}(\phi)=\underset{t\rightarrow0}{\lim}%
{\displaystyle\int\limits_{\text{M}}}
\Upsilon_{t}(\phi,\tau,z)vol(g)_{z}=%
{\displaystyle\int\limits_{\text{M}}}
Tr(\phi)vol(g). \label{As10}%
\end{equation}

\end{theorem}

\textbf{Proof:} Theorem \ref{Trace01} follows directly from Theorem
\ref{Trace00}, Lemma \ref{tr} and the definition of $\Upsilon_{t}(\phi
,\tau,z)$. $\blacksquare$

\subsection{Explicit Formulas}

\begin{theorem}
\label{holder}Let $\mathcal{F}^{^{\prime}}(q,\phi\circ\overline{\psi})$ be
given by the formula $\left(  \ref{v3a}\right)  $. Then for $t>0$ and $q\geq1$
the following equality of the traces of the respective operators
\textit{\ }holds%
\[
Tr\left(  \exp\left(  \left(  -t\triangle_{q-1}^{^{"}}\right)  \circ
\overline{\partial}^{-1}\circ\mathcal{F}^{^{\prime}}(q,\phi\circ\overline
{\psi})\circ\overline{\partial}\right)  \right)  =
\]%
\begin{equation}
Tr\left(  \exp\left(  -t\triangle_{q}^{^{^{\prime}}}\right)  \circ
(\mathcal{F}^{^{\prime}}(q,\phi\circ\overline{\psi}))\right)  . \label{tr10}%
\end{equation}

\end{theorem}

\textbf{Proof: }From Proposition \textbf{2.45} on page 96 in \cite{BGV} it
follows directly that the operators
\[
\exp(-t\Delta_{q-1}^{"})\circ\overline{\partial}^{-1}\circ\mathcal{F}%
^{^{\prime}}(q,\phi\circ\overline{\psi})\circ\overline{\partial},\text{ }%
\exp(-t\Delta_{q}^{^{\prime}})\circ\mathcal{F}^{^{\prime}}(q,\phi
\circ\overline{\psi})
\]
are of trace class since the operators $\exp(-t\Delta_{q-1}^{"})$ have smooth
kernels for $q\geq1$. We know from Proposition \textbf{2.45 }in \cite{BGV}
that we have the following formula:%
\begin{equation}
Tr(DK)=Tr(DA) \label{Com}%
\end{equation}
where $D$ is a differential operator and $A$ is an operator with a smooth
kernel. By using $\left(  \ref{Com}\right)  $ and the fact that the operators
$\Delta_{q}$ and $\overline{\partial}$ commute we derive Theorem \ref{holder}.
Theorem \ref{holder} is proved. $\blacksquare$

\begin{remark}
\label{??}From Definition \ref{Hilb} of the operator $\mathcal{F}^{^{\prime}%
}(q,\phi_{i}\circ\overline{\phi_{j}})$ and Remark \ref{?} we know that it can
be represented pointwise by a matrix which we will denote by $\mathcal{F}%
^{^{\prime}}(q,\left(  \phi_{i}\circ\overline{\phi_{j}}\right)  )$. Since
$k_{t}^{\#}(z,w,0)$ is also a matrix of the same dimension as the operator
$\mathcal{F}^{^{\prime}}(q,\phi_{i}\circ\overline{\phi_{j}})$ we get that the
operator $k_{t}^{\#}(z,w,0)\circ\mathcal{F}^{^{\prime}}(q,\phi_{i}%
\circ\overline{\phi_{j}})$ will be represented pointwise by the product of
finite dimensional matrices. So the integral
\[
\int\limits_{\text{M}}Tr\left(  \left(  k_{t/2}^{\#}(z,w,0)\right)  ^{\ast
}\circ k_{t/2}^{\#}(z,w,0)\circ\mathcal{F}^{^{\prime}}\left(  q,\phi_{i}%
\circ\overline{\phi_{j}}\right)  \right)  vol(g)
\]
makes sense for $t>0$.
\end{remark}

\begin{theorem}
\label{Rtr}Let%
\[
Tr\left(  k_{t}(w,z,\tau)\circ\mathcal{F}^{^{\prime}}(q,\left(  \phi_{i}%
\circ\overline{\phi_{j}}\right)  )\right)  =
\]
\begin{equation}
\frac{\beta_{-n}(\phi_{i}\circ\overline{\phi_{j}})}{t^{n}}+...+\frac
{\beta_{-1}(\phi_{i}\circ\overline{\phi_{j}})}{t}+\beta_{0}(\phi_{i}%
\circ\overline{\phi_{j}})+O(t) \label{ASI}%
\end{equation}
be the short term asymptotic. Then the following limit%
\[
\underset{t\rightarrow0}{\lim}\int\limits_{\text{M}}\left(  \int
\limits_{\text{M}}Tr\left(  \left(  k_{t/2}^{\#}(z,w,0)\right)  ^{\ast}\circ
k_{t/2}^{\#}(z,w,0)\circ\mathcal{F}^{^{\prime}}(q,\phi_{i}\circ\overline
{\phi_{j}}\right)  vol(g)_{w}\right)  vol(g)
\]
exists and%
\[
\beta_{0}(\phi_{i}\circ\overline{\phi_{j}})=
\]%
\[
\underset{t\rightarrow0}{\lim}\int\limits_{\text{M}}\left(  \int
\limits_{\text{M}}Tr\left(  \left(  k_{t/2}^{\#}(z,w,0)\right)  ^{\ast}\circ
k_{t/2}^{\#}(z,w,0)\circ\mathcal{F}^{^{\prime}}(q,\phi_{i}\circ\overline
{\phi_{j}})\right)  vol(g)_{w}\right)  vol(g)=
\]%
\begin{equation}
\int\limits_{\text{M}}Tr\left(  \mathcal{F}^{^{\prime}}\left(  q,\phi_{i}%
\circ\overline{\phi_{j}}\right)  \right)  vol(g)<\infty. \label{Rtr0}%
\end{equation}

\end{theorem}

\textbf{Proof:} Formulas $\left(  \ref{me2a}\right)  $ and $\left(
\ref{me2b}\right)  $ in Theorem \ref{Trace00} imply that%
\[
\beta_{0}(\phi_{i}\circ\overline{\phi_{j}})=\underset{t\rightarrow0}{\lim}%
{\displaystyle\int\limits_{\text{M}}}
\Upsilon_{t}(\mathcal{F}^{^{\prime}}\left(  q,\phi_{i}\circ\overline{\phi_{j}%
}\right)  ,\tau,z)vol(g)_{z}%
\]
Theorem \ref{Trace01} imply that%
\begin{equation}
\beta_{0}(\phi_{i}\circ\overline{\phi_{j}})=\underset{t\rightarrow0}{\lim}%
{\displaystyle\int\limits_{\text{M}}}
\Upsilon_{t}(\mathcal{F}^{^{\prime}}\left(  q,\phi_{i}\circ\overline{\phi_{j}%
}\right)  ,\tau,z)vol(g)_{z}=Tr\left(  \mathcal{F}^{^{\prime}}(q,\phi_{i}%
\circ\overline{\phi_{j}})\right)  |_{z}. \label{Rtr5}%
\end{equation}
Formula $\left(  \ref{Rtr5}\right)  $ implies formula $\left(  \ref{Rtr0}%
\right)  .$ Theorem \ref{Rtr} is proved. $\blacksquare$

\section{The Variational Formulas}

\subsection{Preliminary Formulas}

\begin{lemma}
\label{L1}The following formulas are true for $1\leq q\leq n$:
\begin{equation}
\frac{\partial}{\partial\tau^{i}}\left(  \overline{\partial_{\tau}}\right)
|_{\tau=0}=-F(q,\phi_{i})\circ\partial\label{v3}%
\end{equation}
and
\begin{equation}
\frac{\overline{\partial}}{\overline{\partial\tau^{i}}}\left(  \partial_{\tau
}\right)  |_{\tau=0}=-F(q,\overline{\phi_{i}})\circ\overline{\partial}.
\label{v4}%
\end{equation}

\end{lemma}

\textbf{Proof: }From the expression of $\overline{\partial}_{\tau}$ given in
Definition \ref{tod3}:%

\[
\overline{\partial}_{\tau}=\frac{\overline{\partial}}{\overline{\partial
z^{j}}}-\sum_{m=1}^{N}\left(  \sum_{k=1}^{n}(\phi_{m})_{\overline{j}}^{k}%
\frac{\partial}{\partial z^{k}}\right)  \tau^{m}+O(\tau^{2})),
\]
we conclude that
\begin{equation}
\frac{\partial}{\partial\tau^{i}}\left(  \overline{\partial_{\tau}}\right)
|_{\tau=0}=-\sum_{k=1}^{N}(\phi_{i})_{\overline{j}}^{k}\frac{\partial
}{\partial z^{k}}. \label{w4a}%
\end{equation}
Formula $\left(  \ref{v4}\right)  $ is proved in the same way as formula
$\left(  \ref{w4a}\right)  $. Lemma \ref{L1} follows directly from Definition
\ref{Hilb} of the linear operators F$^{^{\prime}}(q,\phi)$ and F$^{"}%
(q,\phi).$ $\blacksquare$

\begin{corollary}
\label{ip2a}The following formulas are true for $1\leq q\leq n$:
\[
\left(  \frac{\partial}{\partial\tau^{i}}\left(  \triangle_{\tau,q}^{^{"}%
}\right)  |_{\operatorname{Im}\overline{\partial_{\tau}}^{\ast}}\right)
|_{\tau=0}=-\Delta_{0,q}^{"}\circ\overline{\partial_{0}}^{-1}\circ
F(q,\phi_{i})\circ\partial_{0},
\]
and
\begin{equation}
\left(  \frac{\overline{\partial}}{\overline{\partial\tau^{j}}}\left(
\triangle_{\tau,q}^{^{"}}\right)  |_{\operatorname{Im}\overline{\partial
_{\tau}}^{\ast}}\right)  |_{\tau=0}=-\triangle_{0,q}^{^{"}}\circ\partial
_{0}^{-1}\circ F^{^{\prime}}(q+1,\overline{\phi_{i}})\circ\overline
{\partial_{0}}. \label{ah8}%
\end{equation}

\end{corollary}

\textbf{Proof: }From the standard facts of K\"{a}hler geometry we obtain that
on $\operatorname{Im}\overline{\partial}^{\ast}$ in $\Omega_{\text{M}}^{0,q}$
we have
\begin{equation}
\triangle_{\tau,q}^{^{^{"}}}|_{\operatorname{Im}\overline{\partial}^{\ast}%
}=\Lambda_{\tau}\circ\partial_{\tau}\circ\overline{\partial_{\tau}}%
=\Lambda_{\tau}\circ\overline{\partial_{\tau}}\circ\partial_{\tau}.
\label{ah15}%
\end{equation}
We know from $\left(  \ref{v3}\right)  $ and $\left(  \ref{v4}\right)  $ that
\begin{equation}
\frac{\partial}{\partial\tau^{j}}\left(  \overline{\partial_{\tau}}\right)
|_{\tau=0}=-F^{^{\prime}}(q+1,\phi_{j})\circ\partial_{0},\text{ }%
\frac{\partial}{\partial\tau^{j}}\left(  \Lambda_{\tau}\right)  =\frac
{\partial}{\partial\tau^{j}}\left(  \partial_{\tau}\right)  =0. \label{ah12a}%
\end{equation}
Combining $\left(  \ref{ah12a}\right)  $, $\left(  \ref{ah15}\right)  $ we
obtain:%
\[
\frac{\partial}{\partial\tau^{j}}\left(  \triangle_{\tau,q}^{^{"}}\right)
|_{\tau=0}=\frac{\partial}{\partial\tau^{j}}\left(  \Lambda_{\tau}%
\circ\overline{\partial_{\tau}}\circ\overline{\partial_{\tau}}\right)
|_{\tau=0}=
\]%
\begin{equation}
\left(  \overline{\partial_{\tau}}^{\ast}\circ\frac{\partial}{\partial\tau
^{j}}\left(  \overline{\partial_{\tau}}\right)  \right)  |_{\tau=0}%
=-\overline{\partial_{0}}^{\ast}\circ F^{^{\prime}}(q+1,\phi_{j})\circ
\partial_{0}. \label{ah16}%
\end{equation}
Thus on $\operatorname{Im}\overline{\partial}^{\ast}$ we have
\begin{equation}
\partial_{\tau}^{\ast}=\triangle_{\tau,q}^{^{"}}\circ\overline{\partial_{\tau
}}^{-1}. \label{ah17}%
\end{equation}
Substituting $\left(  \ref{ah17}\right)  $ in $\left(  \ref{ah16}\right)  $ we
obtain the first formula in $\left(  \ref{ah8}\right)  $. In the same manner
we obtain the second formula in $\left(  \ref{ah8}\right)  $. Corollary
\ref{ip2a} is proved. $\blacksquare$

\subsection{The Computation of the Antiholomorphic Derivative of $\zeta
_{\tau,q-1}^{"}(s)$}

First we will compute the antiholomorphic derivative of $\zeta_{\tau,q}^{^{"}%
}(s).$

\begin{theorem}
\label{ip}The following formula is true for $t>0$:%
\[
\frac{\overline{\partial}}{\overline{\partial\tau^{i}}}\left(  \zeta_{q,\tau
}^{^{"}}(s)\right)  |_{\tau=0}=
\]
\[
\frac{1}{\Gamma(s)}\int\limits_{0}^{\infty}Tr\left(  \exp(-t(\triangle
_{0,q}^{^{"}})\circ\triangle_{0,q}^{^{"}}\circ\partial_{0}^{-1}\circ
\mathcal{F}^{^{\prime}}(q+1,\overline{\phi_{i}})\circ\overline{\partial_{0}%
}\right)  t^{s}dt.
\]

\end{theorem}

\textbf{Proof: }For the proof of Theorem \ref{ip} we will need the following Lemma:

\begin{lemma}
\label{t11}The following formula is true for $t>0$ and $0<q<n:$%
\[
\frac{\overline{\partial}}{\overline{\partial\tau^{i}}}\left(  Tr(\exp\left(
-t\triangle_{\tau,q}^{^{"}}\right)  \right)  |_{\tau=0}=
\]%
\begin{equation}
tTr\left(  \exp\left(  -t\triangle_{0,q}^{^{"}}\right)  \circ\Delta_{\tau
,q}^{^{"}}\circ\partial^{-1}\circ\mathcal{F}^{^{\prime}}(q+1,\overline
{\phi_{i}})\circ\overline{\partial}\right)  |_{\tau=0}. \label{v7}%
\end{equation}

\end{lemma}

\textbf{Proof: }Direct computations based on Proposition \textbf{9.38.} on
page 304 of the book \cite{DK} show that:%

\begin{equation}
\frac{\overline{\partial}}{\overline{\partial\tau^{i}}}\left(  Tr(\exp\left(
-t\triangle_{\tau,q}^{^{"}}\right)  \right)  |_{\tau=0}=-t\left(  \exp\left(
-t\triangle_{\tau,q}^{^{"}}\right)  \circ\frac{\overline{\partial}}%
{\overline{\partial\tau^{i}}}\left(  \triangle_{\tau,q}^{"}\right)  \right)
|_{\tau=0}. \label{v8}%
\end{equation}
See also \cite{BGV} page 98 Theorem \textbf{2.48. }Formulas $\left(
\ref{v3}\right)  $ and $\left(  \ref{v4}\right)  $ in Lemma \ref{L1} imply
that
\begin{equation}
\frac{\overline{\partial}}{\overline{\partial\tau^{i}}}\left(  \partial_{\tau
}\right)  =-F^{^{\prime}}(q+1,\overline{\phi_{i}(\tau)})\circ\overline
{\partial} \label{v9a}%
\end{equation}
and on $\operatorname{Im}\overline{\partial}^{\ast}$ we have
\begin{equation}
\frac{\overline{\partial}}{\overline{\partial\tau^{i}}}\left(  \partial_{\tau
}^{\ast}\right)  =\frac{\overline{\partial}}{\overline{\partial\tau^{i}}%
}\left(  \Lambda\circ\overline{\partial_{\tau}}\right)  =\left(
\frac{\overline{\partial}}{\overline{\partial\tau^{i}}}\Lambda\right)
\circ\overline{\partial_{\tau}}+\Lambda\circ\frac{\overline{\partial}%
}{\overline{\partial\tau^{i}}}\left(  \overline{\partial_{\tau}}\right)  =0.
\label{v9b}%
\end{equation}
The last equality follows from Corollary. \ref{L1} and \ref{L2}. On K\"{a}hler
manifolds we know that $\partial^{\ast}\circ\partial+\partial\circ
\partial^{\ast}=\overline{\partial}^{\ast}\circ\overline{\partial}%
+\overline{\partial}\circ\overline{\partial}^{\ast}.$ So we deduce that%
\[
\Delta_{\tau,q}^{"}=\left(  \partial_{\tau}^{\ast}\circ\partial_{\tau
}+\partial_{\tau}\circ\partial_{\tau}^{\ast}\right)  |_{\operatorname{Im}%
\overline{\partial_{\tau}}^{\ast}}=\partial_{\tau}^{\ast}\circ\partial_{\tau
}|_{\operatorname{Im}\overline{\partial_{\tau}}^{\ast}}.
\]
Thus from formulas $\left(  \ref{v9a}\right)  $ and $\left(  \ref{v9b}\right)
$ it follows$:$%

\begin{equation}
\frac{\overline{\partial}}{\overline{\partial\tau^{i}}}\left(  \triangle
_{\tau,q}^{^{"}}\right)  =\left(  \partial_{\tau}^{\ast}\circ\frac
{\overline{\partial}}{\overline{\partial\tau^{i}}}\left(  \partial_{\tau
}\right)  \right)  =-\partial_{\tau}^{\ast}\circ F^{^{\prime}}(q+1,\overline
{\phi_{i}(\tau)})\circ\overline{\partial}. \label{v9}%
\end{equation}
By substituting in $\left(  \ref{v8}\right)  $ the expression from $\left(
\ref{v9}\right)  $ we obtain:%

\begin{equation}
\frac{\overline{\partial}}{\overline{\partial\tau^{i}}}\left(  Tr(\exp\left(
-t\triangle_{\tau,q}^{^{"}}\right)  \right)  |_{\tau=0}=tTr\left(  \exp\left(
-t\triangle_{q}^{^{"}}\right)  \circ\partial_{0}^{\ast}\circ F^{^{\prime}%
}(q+1,\overline{\phi_{i}})\circ\overline{\partial_{0}}\right)  . \label{v10}%
\end{equation}
The operator $\partial_{\tau}^{\ast}$ is well defined on the space of
$C^{\infty}$ (0,q) forms on M$_{\tau}.$ So the following formula is true on
$\operatorname{Im}\overline{\partial}_{\tau}^{\ast}$:%

\begin{equation}
\partial_{\tau}^{\ast}=(\triangle_{\tau,q}^{^{"}})\circ\left(  \partial_{\tau
}\right)  ^{-1}. \label{v11}%
\end{equation}
Substituting the expression for $\partial_{\tau}^{\ast}$ in formula $\left(
\ref{v11}\right)  $ in $\left(  \ref{v10}\right)  $, we deduce formula
\ $\left(  \ref{v7}\right)  .$ Lemma \ref{t11} is proved. $\blacksquare$

\textbf{The end of the} \textbf{proof of Theorem }\ref{ip}\textbf{: }The
definition of the zeta function implies that%
\[
\frac{\overline{\partial}}{\overline{\partial\tau^{i}}}\left(  \zeta
_{\Delta_{\tau,q}^{"}}(s)\right)  |_{\tau=0}=\frac{\overline{\partial}%
}{\overline{\partial\tau^{i}}}\left(  \frac{1}{\Gamma(s)}\int\limits_{0}%
^{\infty}Tr\left(  \exp(-t(\triangle_{\tau,q}^{^{^{"}}})\right)
t^{s-1}dt\right)  |_{\tau=0}=
\]%
\begin{equation}
\frac{1}{\Gamma(s)}\int\limits_{0}^{\infty}\left(  \frac{\overline{\partial}%
}{\overline{\partial\tau^{i}}}\left(  Tr\exp(-t(\triangle_{\tau,q}^{^{^{"}}%
})\right)  \right)  |_{\tau=0}t^{s-1}dt. \label{ah0}%
\end{equation}
Substituting in $\left(  \ref{ah0}\right)  $ the expression for $\frac
{\overline{\partial}}{\overline{\partial\tau^{i}}}\left(  Tr\left(
\exp(-t(\triangle_{\tau,q}^{"})\right)  \right)  $ given by $\left(
\ref{v7}\right)  $ we obtain:%

\[
\frac{\overline{\partial}}{\overline{\partial\tau^{i}}}\left(  \zeta
_{\Delta_{\tau,q}^{"}}(s)\right)  |_{\tau=0}=
\]%
\begin{equation}
\frac{1}{\Gamma(s)}\int\limits_{0}^{\infty}Tr\left(  \left(  \exp
(-t(\triangle_{0,q}^{"})\right)  \circ\left(  \triangle_{0,q}^{"}\right)
\circ\partial^{-1}\circ F^{^{\prime}}(q+1,\overline{\phi_{i}})\circ
\overline{\partial}\right)  t^{s}dt. \label{ah-1}%
\end{equation}
Theorem \ref{ip} is proved. $\blacksquare$

\subsection{The Computation of the Hessian of $\zeta_{\tau,q}^{^{"}}(s)$}

\begin{theorem}
\label{ip2}The following formula holds:%
\[
\frac{\partial^{2}}{\partial\tau^{j}\overline{\partial\tau^{i}}}\left(
\zeta_{\tau,q}^{^{"}}(s)\right)  |_{\tau=0}=
\]
\begin{equation}
\frac{\left(  1-s\right)  s}{\Gamma(s)}\int\limits_{0}^{\infty}Tr\left(
\left(  \exp(-t(\triangle_{0,q+1}^{^{^{^{\prime}}}})\right)  \circ
\mathcal{F}^{^{\prime}}(q+1,\phi_{j}\circ\overline{\phi_{i}})\right)
t^{s-1}dt. \label{ah10}%
\end{equation}

\end{theorem}

\textbf{Proof: }The facts that the operators
\[
\frac{\partial}{\partial\tau^{i}}\left(  \overline{\partial_{\tau}}\right)
=\left(  -\phi_{i}+O(\tau)\right)  \circ\partial_{0}%
\]
depend holomorphically on $\tau$ and the operator $\partial_{\tau}^{-1}$
depends antiholomorphically imply that the operators $\partial_{\tau}%
^{-1}\circ F^{^{\prime}}(q+1,\overline{\phi_{i}})\circ\overline{\partial_{0}}$
depend antiholomorphically on the coordinates $\tau=(\tau^{1},...,\tau^{N}).$
By using the explicit formula $\left(  \ref{ah-1}\right)  $ for the
antiholomorphic derivative of $\zeta_{\tau,q}^{"}(s)$ and that%
\[
\frac{\partial}{\partial\tau^{j}}\left(  \partial_{\tau}^{-1}\circ
F^{^{\prime}}(q+1,\overline{\phi_{i}})\circ\overline{\partial_{0}}\right)  =0
\]
we derive
\[
\frac{\partial^{2}}{\overline{\partial\tau^{j}}\partial\tau^{i}}\left(
\zeta_{\tau,q}^{"}(s)\right)  =
\]%
\[
\frac{1}{\Gamma(s)}\int\limits_{0}^{\infty}Tr\left(  \left(  \frac{\partial
}{\partial\tau^{j}}\exp(-t\left(  \triangle_{\tau,q}^{^{"}}\right)  \right)
\circ\left(  \triangle_{\tau,q}^{^{"}}\right)  \circ\partial_{\tau}^{-1}\circ
F^{^{\prime}}(q+1,\overline{\phi_{i}})\circ\overline{\partial_{0}}\right)
t^{s}dt+
\]%
\begin{equation}
\frac{1}{\Gamma(s)}\int\limits_{0}^{\infty}Tr\left(  \left(  \exp
(-t(\triangle_{\tau,q}^{"})\right)  \circ\frac{\partial}{\partial\tau^{j}%
}\left(  \triangle_{\tau,q}^{^{"}}\right)  \circ\partial_{0}^{-1}\circ
F^{^{\prime}}(q+1,\overline{\phi_{i}})\circ\overline{\partial_{0}}\right)
t^{s}dt. \label{ah4}%
\end{equation}

\begin{lemma}
\label{AH1}We have the following expression:%
\[
\frac{1}{\Gamma(s)}\int\limits_{0}^{\infty}Tr\left(  \left(  \frac{\partial
}{\partial\tau^{j}}\left(  \exp(-t\left(  \triangle_{\tau,q}^{^{"}}\right)
\right)  \right)  \circ\left(  \triangle_{\tau,q}^{^{"}}\right)  \circ
\partial_{\tau}^{-1}\circ F^{^{\prime}}(q+1,\overline{\phi_{i}})\circ
\overline{\partial_{0}}\right)  t^{s}dt=
\]%
\begin{equation}
\frac{-s}{\Gamma(s)}\int\limits_{0}^{\infty}Tr\left(  \left(  \exp
(-t\triangle_{0,q}^{^{"}}\right)  \circ\frac{\partial}{\partial\tau^{j}%
}\left(  \triangle_{\tau,q}^{^{"}}\right)  \circ\partial_{\tau}^{-1}\circ
F^{^{\prime}}(q+1,\overline{\phi_{i}})\circ\overline{\partial_{0}}\right)
t^{s}dt. \label{AH10}%
\end{equation}

\end{lemma}

\textbf{Proof: }Direct computations show that
\[
\frac{1}{\Gamma(s)}\left(  \int\limits_{0}^{\infty}Tr\left(  \left(
\frac{\partial}{\partial\tau^{j}}\left(  \exp(-t\left(  \triangle_{\tau
,q}^{^{"}}\right)  \right)  \right)  \circ\left(  \triangle_{\tau,q}^{^{"}%
}\right)  \circ\partial_{\tau}^{-1}\circ F^{^{\prime}}(q+1,\overline{\phi_{i}%
})\circ\overline{\partial_{0}}\right)  t^{s}dt\right)  |_{\tau=0}=
\]%
\begin{equation}
\frac{-1}{\Gamma(s)}\int\limits_{0}^{\infty}\left(  \frac{d}{dt}Tr\left(
\left(  \frac{\partial}{\partial\tau^{j}}\left(  \exp(-t\left(  \triangle
_{\tau,q}^{^{"}}\right)  \right)  \right)  |_{\tau=0}\circ\partial_{0}%
^{-1}\circ F^{^{\prime}}(q+1,\overline{\phi_{i}})\circ\overline{\partial_{0}%
}\right)  \right)  t^{s}dt \label{ah4a}%
\end{equation}
By integrating by parts the right hand side of formula $\left(  \ref{ah4a}%
\right)  $ we deduce that:
\[
\frac{-1}{\Gamma(s)}\int\limits_{0}^{\infty}\left(  \frac{d}{dt}Tr\left(
\left(  \frac{\partial}{\partial\tau^{j}}\left(  \exp(-t\left(  \triangle
_{\tau,q}^{^{"}}\right)  \right)  \right)  |_{\tau=0}\circ\partial_{0}%
^{-1}\circ F^{^{\prime}}(q+1,\overline{\phi_{i}})\circ\overline{\partial_{0}%
}\right)  \right)  t^{s}dt=
\]%
\begin{equation}
\frac{s}{\Gamma(s)}\int\limits_{0}^{\infty}Tr\left(  \left(  \frac{\partial
}{\partial\tau^{j}}\left(  \exp(-t\left(  \triangle_{\tau,q}^{^{"}}\right)
\right)  \right)  |_{\tau=0}\circ\partial_{0}^{-1}\circ F^{^{\prime}%
}(q+1,\overline{\phi_{i}})\circ\overline{\partial_{0}}\right)  t^{s-1}dt.
\label{ah4b}%
\end{equation}
Direct computations of the right hand side of $\left(  \ref{ah4b}\right)  $
show that:
\[
\frac{s}{\Gamma(s)}\int\limits_{0}^{\infty}Tr\left(  \left(  \frac{\partial
}{\partial\tau^{j}}\left(  \exp(-t\triangle_{\tau,q}^{^{"}}\right)  \right)
|_{\tau=0}\circ\partial_{0}^{-1}\circ F^{^{\prime}}(q+1,\overline{\phi_{i}%
})\circ\overline{\partial_{0}}\right)  t^{s-1}dt=
\]%
\[
\frac{-s}{\Gamma(s)}\int\limits_{0}^{\infty}Tr\left(  \left(  \exp
(-t\triangle_{\tau,q}^{^{"}}\right)  \circ\frac{\partial}{\partial\tau^{j}%
}\left(  \triangle_{\tau,q}^{^{"}}\right)  |_{\tau=0}\circ\partial_{0}%
^{-1}\circ F^{^{\prime}}(q+1,\overline{\phi_{i}})\circ\overline{\partial_{0}%
}\right)  t^{s}dt.
\]
Formula $\left(  \ref{AH10}\right)  $ is proved. $\blacksquare$

Substituting in $\left(  \ref{ah4}\right)  $ the expression from $\left(
\ref{AH10}\right)  $ for
\[
\frac{1}{\Gamma(s)}\int\limits_{0}^{\infty}Tr\left(  \left(  \frac{\partial
}{\partial\tau^{j}}\left(  \exp(-t\left(  \triangle_{\tau,q}^{^{"}}\right)
\right)  \right)  \circ\left(  \triangle_{\tau,q}^{^{"}}\right)  \circ
\partial_{\tau}^{-1}\circ F^{^{\prime}}(q+1,\overline{\phi_{i}})\circ
\overline{\partial_{0}}\right)  t^{s}dt
\]
we get:
\[
\frac{\partial^{2}}{\overline{\partial\tau^{j}}\partial\tau^{i}}\left(
\zeta_{q-1,\tau}^{"}(s)\right)  |_{\tau=0}=
\]%
\begin{equation}
\frac{1-s}{\Gamma(s)}\int\limits_{0}^{\infty}Tr\left(  \left(  \exp(-t\left(
\triangle_{\tau,q}^{^{"}}\right)  \right)  \circ\frac{\partial}{\partial
\tau^{j}}\left(  \triangle_{\tau,q}^{^{"}}\right)  |_{\tau=0}\circ\partial
_{0}^{-1}\circ F^{^{\prime}}(q+1,\overline{\phi_{i}})\circ\overline
{\partial_{0}}\right)  t^{s}dt. \label{ah40}%
\end{equation}

\begin{lemma}
\label{AH2}%
\[
\int\limits_{0}^{\infty}Tr\left(  \left(  \exp(-t\left(  \triangle_{\tau
,q}^{^{"}}\right)  \right)  \circ\frac{\partial}{\partial\tau^{j}}\left(
\triangle_{\tau,q}^{^{"}}\right)  |_{\tau=0}\circ\partial_{0}^{-1}\circ
F^{^{\prime}}(q+1,\overline{\phi_{i}})\circ\overline{\partial_{0}}\right)
t^{s}dt=
\]%
\begin{equation}
s\int\limits_{0}^{\infty}Tr\left(  \left(  \exp(-t\left(  \triangle_{\tau
,q}^{^{"}}\right)  \right)  |_{\tau=0}\circ\left(  \overline{\partial_{0}%
}\right)  ^{-1}\circ\mathcal{F}^{^{\prime}}(q+1,\phi_{j}\circ\overline
{\phi_{i}})\circ\overline{\partial_{0}}\right)  t^{s-1}dt. \label{AH20}%
\end{equation}

\end{lemma}

\textbf{Proof: }Substituting the expression of $\left(  \ref{ah8}\right)  $
for
\[
\frac{\partial}{\partial\tau^{j}}\left(  \triangle_{\tau,q}^{^{"}}\right)
|_{\tau=0}=-\triangle_{0,q}^{^{"}}\circ\left(  \overline{\partial_{0}}\right)
^{-1}\circ F^{^{\prime}}(q+1,\phi_{j})\circ\overline{\partial_{0}}%
\]
in the expression of
\[
\int\limits_{0}^{\infty}Tr\left(  \left(  \exp(-t\left(  \triangle_{\tau
,q}^{^{"}}\right)  \right)  \circ\frac{\partial}{\partial\tau^{j}}\left(
\triangle_{\tau,q}^{^{"}}\right)  |_{\tau=0}\circ\partial_{0}^{-1}\circ
F^{^{\prime}}(q+1,\overline{\phi_{i}})\circ\overline{\partial_{0}}\right)
t^{s}dt
\]
we obtain:%
\[
\int\limits_{0}^{\infty}Tr\left(  \left(  \exp(-t\left(  \triangle_{\tau
,q}^{^{"}}\right)  \right)  \circ\frac{\partial}{\partial\tau^{j}}\left(
\triangle_{\tau,q}^{^{"}}\right)  |_{\tau=0}\circ\partial_{0}^{-1}\circ
F^{^{\prime}}(q+1,\overline{\phi_{i}})\circ\overline{\partial_{0}}\right)
t^{s}dt=
\]%
\begin{equation}
-\int\limits_{0}^{\infty}Tr\left(  \left(  \exp(-t\left(  \triangle
_{0,q}^{^{"}}\right)  \right)  \circ\triangle_{0,q}^{^{"}}\circ\partial
_{0}^{-1}\circ F^{^{\prime}}(q+1,\overline{\phi_{i}})\circ F^{^{\prime}%
}(q+1,\overline{\phi_{i}})\circ\overline{\partial_{0}}\right)  t^{s}dt.
\label{AH21a}%
\end{equation}
Simple observations show that
\[
\int\limits_{0}^{\infty}Tr\left(  \left(  \exp(-t\left(  \triangle_{0,q}%
^{^{"}}\right)  \right)  \circ\triangle_{0,q}^{^{"}}\circ\partial_{0}%
^{-1}\circ F^{^{\prime}}(q+1,\overline{\phi_{i}})\circ F^{^{\prime}%
}(q+1,\overline{\phi_{i}})\circ\overline{\partial_{0}}\right)  t^{s}dt=
\]%
\begin{equation}
-\int\limits_{0}^{\infty}Tr\left(  \frac{d}{dt}\exp(-t\left(  \triangle
_{0,q}^{^{"}}\right)  \partial_{0}^{-1}\circ F^{^{\prime}}(q+1,\overline
{\phi_{i}})\circ F^{^{\prime}}(q+1,\overline{\phi_{i}})\circ\overline
{\partial_{0}}\right)  t^{s}dt. \label{AH21b}%
\end{equation}
By integrating by parts the right hand side of $\left(  \ref{AH21b}\right)  $
we obtain:%
\[
-\int\limits_{0}^{\infty}Tr\left(  \frac{d}{dt}\exp(-t\left(  \triangle
_{0,q}^{^{"}}\right)  \partial_{0}^{-1}\circ F^{^{\prime}}(q+1,\overline
{\phi_{i}})\circ F^{^{\prime}}(q+1,\overline{\phi_{i}})\circ\overline
{\partial_{0}}\right)  t^{s}dt=
\]%
\[
s\int\limits_{0}^{\infty}Tr\left(  \exp(-t\left(  \triangle_{0,q}^{^{"}%
}\right)  \partial_{0}^{-1}\circ F^{^{\prime}}(q+1,\overline{\phi_{i}})\circ
F^{^{\prime}}(q+1,\overline{\phi_{i}})\circ\overline{\partial_{0}}\right)
t^{s-1}dt.
\]
Thus we derive formula $\left(  \ref{AH20}\right)  ,$ i.e.%
\[
\int\limits_{0}^{\infty}Tr\left(  \left(  \exp(-t\left(  \triangle_{\tau
,q}^{^{"}}\right)  \right)  \circ\frac{\partial}{\partial\tau^{j}}\left(
\triangle_{\tau,q}^{^{"}}\right)  |_{\tau=0}\circ\partial_{0}^{-1}\circ
F^{^{\prime}}(q+1,\overline{\phi_{i}})\circ\overline{\partial_{0}}\right)
t^{s}dt=
\]%
\[
s\int\limits_{0}^{\infty}Tr\left(  \left(  \exp(-t\left(  \triangle_{\tau
,q}^{^{"}}\right)  \right)  |_{\tau=0}\circ\left(  \overline{\partial_{0}%
}\right)  ^{-1}\circ\mathcal{F}^{^{\prime}}(q+1,\phi_{j}\circ\overline
{\phi_{i}})\circ\overline{\partial_{0}}\right)  t^{s-1}dt.
\]
Lemma \ref{AH2} is proved. $\blacksquare$

Substituting the expression of $\left(  \ref{AH20}\right)  $ in the expression
$\left(  \ref{ah40}\right)  $ we get the following equality:
\[
\frac{\partial^{2}}{\overline{\partial\tau^{j}}\partial\tau^{i}}\left(
\zeta_{q,\tau}^{"}(s)\right)  |_{\tau=0}=
\]%
\begin{equation}
\frac{\left(  1-s\right)  s}{\Gamma(s)}\int\limits_{0}^{\infty}Tr\left(
\left(  \exp(-t(\triangle_{0,q}^{^{"}})\right)  \circ\overline{\partial_{0}%
}^{-1}\circ\mathcal{F}^{^{\prime}}(q+1,\phi_{j}\circ\overline{\phi_{i}}%
)\circ\overline{\partial_{0}}\right)  t^{s-1}dt. \label{ah20}%
\end{equation}
Applying Theorem \ref{holder} we deduce Theorem \ref{ip2}. $\blacksquare$

\subsection{The Computations of the Hessian of $\log\det\Delta_{\tau,q}$}

\begin{theorem}
\label{Main}The following formula is true:
\[
\frac{\partial^{2}}{\overline{\partial\tau^{j}}\partial\tau^{i}}\left(
\log\det\Delta_{\tau,q}^{"}\right)  |_{\tau=0}=
\]%
\[
-\underset{t\rightarrow0}{\lim}\int\limits_{\text{M}}\left(  \int
\limits_{\text{M}}Tr\left(  4\pi t)^{-\frac{n}{2}}\exp\left(  -\frac{d_{0}%
^{2}(w,z)}{4t}\right)  \circ\mathcal{P}\circ\mathcal{F}^{^{\prime}}\left(
q+1,\phi_{i}\circ\overline{\phi_{j}}\right)  \right)  vol(g)_{w}\right)
vol(g)_{z}=
\]%
\begin{equation}
-\int\limits_{\text{M}}Tr\left(  \mathcal{F}^{^{\prime}}\left(  q+1,\left(
\phi_{i}\circ\overline{\phi_{j}}\right)  \right)  )\right)  vol(g).
\label{Rtr7}%
\end{equation}

\end{theorem}

\textbf{Proof: }$\zeta_{\tau,q}^{"}(s)$ is obtained from the meromorphic
continuation of
\[
\frac{1}{\Gamma(s)}%
{\displaystyle\int\limits_{0}^{\infty}}
Tr\left(  \exp(-t\Delta_{\tau,q}^{"})\right)  t^{s-1}dt.
\]
Thus it is a meromorphic function on $\mathbb{C}$ well defined at $0.$ So we
get that $\zeta_{\tau,q}^{"}(s)=\mu_{0}(\tau)+\mu_{1}(\tau)s+O(s^{2}).$ From
here we deduce%
\[
\frac{\partial^{2}}{\overline{\partial\tau^{j}}\partial\tau^{i}}\left(
\zeta_{\tau,q}^{"}(s)\right)  |_{\tau=0}=\frac{\partial^{2}}{\overline
{\partial\tau^{j}}\partial\tau^{i}}\mu_{0}(\tau)|_{\tau=0}+\left(
\frac{\partial^{2}}{\overline{\partial\tau^{j}}\partial\tau^{i}}\mu_{1}%
(\tau)\right)  |_{\tau=0}s+O(s^{2})=
\]%
\begin{equation}
\alpha_{0}+\alpha_{1}s+O(s^{2}). \label{Rtr12}%
\end{equation}
Thus from the definition of the regularized determinant
\[
\log\det(\triangle_{\tau,q}^{"})=\left(  \frac{d}{ds}\left(  -\zeta_{\tau
,q}^{"}(s)\right)  \right)  |_{s=0}%
\]
we see that
\begin{equation}
\frac{\partial^{2}}{\overline{\partial\tau^{j}}\partial\tau^{i}}\left(
\log\det\Delta_{\tau,q}^{"}\right)  |_{\tau=0}=\frac{d}{ds}\left(
\frac{\partial^{2}}{\overline{\partial\tau^{j}}\partial\tau^{i}}\left(
\left(  -\zeta_{\tau,q}^{"}(s)\right)  \right)  |_{\tau=0}\right)
|_{s=0}=-\alpha_{1}. \label{Rtr11}%
\end{equation}
Combining formula $\left(  \ref{ah10}\right)  $
\[
\frac{\partial^{2}}{\partial\tau^{j}\overline{\partial\tau^{i}}}\left(
\zeta_{\tau,q}^{^{"}}(s)\right)  |_{\tau=0}=\frac{s}{\Gamma(s)}\int
\limits_{0}^{\infty}Tr\left(  \left(  \exp(-t(\triangle_{0,q+1}^{^{\prime}%
})\right)  \circ\mathcal{F}^{^{\prime}}(q+1,\phi_{j}\circ\overline{\phi_{i}%
})\right)  t^{s-1}dt+
\]%
\[
-\frac{s^{2}}{\Gamma(s)}\int\limits_{0}^{\infty}Tr\left(  \left(
\exp(-t(\triangle_{0,q+1}^{^{^{^{\prime}}}})\right)  \circ\mathcal{F}%
^{^{\prime}}(q+1,\phi_{j}\circ\overline{\phi_{i}})\right)  t^{s-1}dt
\]
with the short term expansion:%
\begin{equation}
Tr\left(  \left(  \exp(-t(\triangle_{0,q+1}^{^{^{\prime}}})\right)
\circ\mathcal{F}^{^{\prime}}(q+1,\phi_{j}\circ\overline{\phi_{i}})\right)  =%
{\displaystyle\sum\limits_{k=-n}^{1}}
\frac{\nu_{k}}{t^{k}}+\nu_{0}+\psi(t) \label{As11}%
\end{equation}
where
\[
\psi(t)=Tr\left(  \left(  \exp(-t(\triangle_{0,q+1}^{^{^{\prime}}})\right)
\circ\mathcal{F}^{^{\prime}}(q+1,\phi_{j}\circ\overline{\phi_{i}})\right)  -%
{\displaystyle\sum\limits_{k=-n}^{1}}
\frac{\nu_{k}}{t^{k}}+\nu_{0}%
\]
we obtain:%
\[
\frac{\partial^{2}}{\partial\tau^{j}\overline{\partial\tau^{i}}}\left(
\zeta_{\tau,q}^{^{"}}(s)\right)  |_{\tau=0}=
\]%
\[
\frac{s}{\Gamma(s)}\left(
{\displaystyle\int\limits_{0}^{1}}
\left(
{\displaystyle\sum\limits_{k=-n}^{1}}
\frac{\nu_{k}}{t^{k}}\right)  t^{s-1}dt+\nu_{0}%
{\displaystyle\int\limits_{0}^{1}}
t^{s-1}dt+%
{\displaystyle\int\limits_{0}^{1}}
\psi(t)t^{s-1}dt\right)  +
\]%
\[
\frac{s}{\Gamma(s)}\left(
{\displaystyle\int\limits_{1}^{\infty}}
Tr\left(  \left(  \exp(-t(\triangle_{0,q+1}^{^{^{\prime}}})\right)
\circ\mathcal{F}^{^{\prime}}(q+1,\phi_{j}\circ\overline{\phi_{i}})\right)
t^{s-1}dt\right)  -
\]%
\[
-\frac{s^{2}}{\Gamma(s)}\left(
{\displaystyle\int\limits_{0}^{1}}
\left(
{\displaystyle\sum\limits_{k=-n}^{1}}
\frac{\nu_{k}}{t^{k}}\right)  t^{s-1}dt+\nu_{0}%
{\displaystyle\int\limits_{0}^{1}}
t^{s-1}dt+%
{\displaystyle\int\limits_{0}^{1}}
\psi(t)t^{s-1}dt\right)  -
\]%
\begin{equation}
-\frac{s^{2}}{\Gamma(s)}\left(
{\displaystyle\int\limits_{1}^{\infty}}
Tr\left(  \left(  \exp(-t(\triangle_{0,q+1}^{^{^{\prime}}})\right)
\circ\mathcal{F}^{^{\prime}}(q+1,\phi_{j}\circ\overline{\phi_{i}})\right)
t^{s-1}dt\right)  . \label{AS11}%
\end{equation}
By using formula $\left(  \ref{AS11}\right)  $ we will prove the following Lemma:

\begin{lemma}
\label{AS00}We have the following formula:%
\[
\frac{\partial^{2}}{\overline{\partial\tau^{j}}\partial\tau^{i}}\left(
\log\det\Delta_{\tau,q}^{"}\right)  |_{\tau=0}=-%
{\displaystyle\int\limits_{\text{M}}}
Tr\mathcal{F}^{^{\prime}}(q+1,\phi_{i}\circ\overline{\phi_{j}})vol(g)=-\alpha
_{1}.
\]

\end{lemma}

\textbf{Proof: }Lemma \textbf{9.34} on page 300 of \cite{BGV} or direct
computations show that for $|s|<\varepsilon$ we have the following identity:%
\[
\frac{1}{\Gamma(s)}\int\limits_{0}^{\infty}Tr\left(  \left(  \exp
(-t(\triangle_{0,q}^{^{\prime}})\right)  \circ\mathcal{F}^{^{\prime}}%
(q+1,\phi_{j}\circ\overline{\phi_{i}})\right)  t^{s-1}dt=
\]%
\[
\frac{1}{\Gamma(s)}\left(
{\displaystyle\int\limits_{0}^{1}}
\left(
{\displaystyle\sum\limits_{k=-n}^{1}}
\frac{\nu_{k}}{t^{k}}\right)  t^{s-1}dt+\nu_{0}%
{\displaystyle\int\limits_{0}^{1}}
t^{s-1}dt+%
{\displaystyle\int\limits_{0}^{1}}
\psi(t)t^{s-1}dt\right)  +
\]%
\begin{equation}
+\frac{1}{\Gamma(s)}%
{\displaystyle\int\limits_{1}^{\infty}}
Tr\left(  \left(  \exp(-t(\triangle_{0,q}^{^{^{\prime}}})\right)
\circ\mathcal{F}^{^{\prime}}(q+1,\phi_{j}\circ\overline{\phi_{i}})\right)
t^{s-1}dt=\frac{\nu_{0}}{s}+\kappa+O(s). \label{AS11a}%
\end{equation}
Combining the expression in $\left(  \ref{AS11a}\right)  $ with the following
standard fact $\frac{s}{\Gamma(s)}=s^{2}+O(s^{3})$ we obtain from formulas
$\left(  \ref{AS11}\right)  $ and $\left(  \ref{AS11a}\right)  $\ that for
$|s|<\varepsilon$%
\begin{equation}
\frac{\partial^{2}}{\partial\tau^{j}\overline{\partial\tau^{i}}}\left(
\zeta_{\tau,q}^{^{"}}(s)\right)  |_{\tau=0}=\nu_{0}s+O(s^{2}). \label{AS11b}%
\end{equation}
Thus according to $\left(  \ref{Rtr11}\right)  $ and $\left(  \ref{AS11b}%
\right)  $
\begin{equation}
\frac{d}{ds}\left(  \frac{\partial^{2}}{\partial\tau^{j}\overline{\partial
\tau^{i}}}\left(  \zeta_{\tau,q}^{^{"}}(s)\right)  |_{\tau=0}\right)
|_{s=0}=\nu_{0}=\alpha_{1}. \label{AS12}%
\end{equation}
Applying Theorem \ref{Trace01} to formula $\left(  \ref{AS12}\right)  $ we
deduce that
\[
\alpha_{1}=\nu_{0}=\int\limits_{\text{M}}Tr\left(  \mathcal{F}^{^{\prime}%
}\left(  q+1,\phi_{i}\circ\overline{\phi_{j}}\right)  \right)  vol(g).
\]
Lemma \ref{AS00} is proved. $\blacksquare$

Lemma \ref{AS00}\ implies directly Theorem \ref{Main}. $\blacksquare$

\subsection{Some Applications of the Variational Formulas}

\begin{theorem}
\label{IZS1}The following identity holds $dd^{c}\left(  \log\det\Delta
_{\tau,1}\right)  =dd^{c}\left(  \log\det\Delta_{\tau,1}^{^{\prime}}\det
\Delta_{\tau,1}^{"}\right)  =-\operatorname{Im}W.P.$
\end{theorem}

\textbf{Proof: }The proof of Theorem \ref{IZS1} is based on the following
formulas which holds for K\"{a}hler manifolds:%
\begin{equation}
\Delta_{\partial}=\Delta\overline{_{\partial}},\text{ }\partial^{\ast}%
=-\ast\overline{\partial}\ast\text{ and }\overline{\partial}^{\ast}%
=-\ast\partial\ast\text{,} \label{du1}%
\end{equation}
where $\ast$ is the Hodge star operator. See \cite{KM} page 95. On CY manifold
we have the following duality: $\ast:\Omega_{\text{M}}^{0,q}\approxeq
\Omega_{\text{M}}^{0,n-q}$ induced by the Hodge star operator $\ast$ of a CY
metric and the holomorphic $n$ form. Using this duality direct check shows
that on CY manifolds we have%
\begin{equation}
\ast\left(  \operatorname{Im}\overline{\partial_{q}}\right)
=\operatorname{Im}\left(  \overline{\partial_{n-q}}^{\ast}\right)  \text{ and
}\ast\left(  \operatorname{Im}\overline{\partial_{q}}^{\ast}\right)
=\operatorname{Im}\left(  \overline{\partial_{n-q}}\right)  \label{du}%
\end{equation}
Formulas $\left(  \ref{du1}\right)  $ and $\left(  \ref{du}\right)  $ imply
that we have $\det\Delta_{\overline{\partial},q}^{^{\prime}}=\det
\Delta_{\partial,q}^{^{\prime}}=\det\Delta_{\partial,n-q}^{"^{\prime}}%
=\det\Delta_{\overline{\partial},n-q}^{"^{\prime}}.$

\begin{lemma}
\label{IZS}We have the following relations between the operators on a CY
manifold%
\[
Tr\left(  \exp(-t\Delta_{\tau,1}^{^{\prime}})\circ\mathcal{F}^{^{\prime}%
}\left(  1,\phi_{i}\circ\overline{\phi_{j}}\right)  \right)  =Tr\left(
\exp(-t\Delta_{\tau,1}^{"})\circ\mathcal{F}^{^{"}}\left(  n-1,\phi_{i}%
\circ\overline{\phi_{j}}\right)  \right)
\]
and
\begin{equation}
\text{ }Tr\left(  \exp(-t\Delta_{\tau,1}^{^{"}})\circ\mathcal{F}^{"}\left(
1,\phi_{i}\circ\overline{\phi_{j}}\right)  \right)  =Tr\left(  \exp
(-t\Delta_{\tau,1}^{^{\prime}})\circ\mathcal{F}^{^{\prime}}\left(
n-1,\phi_{i}\circ\overline{\phi_{j}}\right)  \right)  \label{iz3}%
\end{equation}
by identifying the Hilbert spaces $\operatorname{Im}\overline{\partial}\subset
L^{2}($M$,\Omega_{\text{M}}^{0,1})$ and $\operatorname{Im}\overline{\partial
}^{\ast}\subset L^{2}\left(  \text{M},\Omega_{\text{M}}^{0,1}\right)  $ with
$\operatorname{Im}\partial^{\ast}\subset L^{2}\left(  \text{M},\Omega
_{\text{M}}^{0,n-1}\right)  $ and $\operatorname{Im}\partial\subset
L^{2}\left(  \text{M},\Omega_{\text{M}}^{0,n-1}\right)  $ respectively by
using $\left(  \ref{du}\right)  .$
\end{lemma}

\textbf{Proof:} We will need the following Propositions to prove Lemma
\ref{IZS}:

\begin{proposition}
\label{IZS0}Let $\left\{  \omega_{i}\right\}  $ be an orthonormal basis at
$\Omega_{x}^{1,0}.$ Let $\overline{\omega}_{n}=\overline{\omega_{1}}%
\wedge...\wedge\overline{\omega_{n}}$ be the antiholomorphic volume form. Then
for $\alpha\neq\beta$ we have%
\begin{equation}
\left(  \phi_{i}\circ\overline{\phi_{j}}\wedge id_{n-2}\right)  \left(
\overline{\omega}_{n}\lrcorner\frac{\overline{\partial}}{\overline{\partial
z^{\alpha}}}\right)  =%
{\displaystyle\sum\limits_{\beta}}
\left(  \phi_{i}\circ\overline{\phi_{j}}\right)  _{\overline{\beta}%
}^{\overline{\alpha}}\left(  \overline{\omega}_{n}\lrcorner\frac
{\overline{\partial}}{\overline{\partial z^{\alpha}}}\wedge\frac
{\overline{\partial}}{\overline{\partial z^{\beta}}}\right)  \wedge
\overline{dz^{\alpha}}.\label{L0}%
\end{equation}

\end{proposition}

\textbf{Proof:} Formula \ref{L0} follows directly from the definition of the
linear operator $\left(  \phi_{i}\circ\overline{\phi_{j}}\wedge id_{n-2}%
\right)  .\blacksquare$

\begin{proposition}
\label{IZS2}Let us define $\left(  \phi_{i}\circ\overline{\phi_{j}}\wedge
id_{n-2}\right)  ^{\ast}:\Omega_{\text{M}}^{0,n-1}\rightarrow\Omega_{\text{M}%
}^{0,n-1}$ as follows:%
\begin{equation}
\left(  \phi_{i}\circ\overline{\phi_{j}}\wedge id_{n-2}\right)  ^{\ast}\left(
\ast\overline{\omega_{k}}\right)  =\ast\left(  \left(  \phi_{i}\circ
\overline{\phi_{j}}\right)  \left(  \overline{\omega_{k}}\right)  \right)
\label{LO2}%
\end{equation}
Let us denote by $\emph{M}\left(  \overline{\phi_{j}}\circ\phi_{i}\right)  $
and $\emph{M}\left(  \left(  \phi_{i}\circ\overline{\phi_{j}}\wedge
id_{n-2}\right)  ^{\ast}\right)  $ the matrices of the operators
$\overline{\phi_{j}}\circ\phi_{i}$ and $\left(  \phi_{i}\circ\overline
{\phi_{j}}\wedge id_{n-2}\right)  ^{\ast}$ in the orthonormal bases with
respect to the CY metric.\textbf{ }Then fibrewise we have the equality
\begin{equation}
\emph{M}\left(  \overline{\phi_{j}}\circ\phi_{i}\right)  =\emph{M}\left(
\left(  \phi_{i}\circ\overline{\phi_{j}}\wedge id_{n-2}\right)  ^{\ast
}\right)  . \label{LO1}%
\end{equation}

\end{proposition}

\textbf{Proof:} We need to compute the matrix of the operator $\left(
\phi_{i}\circ\overline{\phi_{j}}\wedge id_{n-2}\right)  ^{\ast}$ in the
orthonormal basis $\overline{\omega_{i_{1}}}\wedge...\wedge\overline
{\omega_{i_{n-k}}}$ and compare it with the matrix of the operator
$\overline{\phi_{j}}\circ\phi_{i}~$of the bundle $\Omega_{\text{M}}^{1,0}$
written in the orthonormal basis $\left\{  \omega_{i}\right\}  .$

Let $\left\{  \omega_{i}\right\}  $ be an orthonormal basis at $\Omega
_{x}^{1,0}.$ According to Lemma \ref{sym} the operators $\phi_{i}:\Omega
_{x}^{1,0}\rightarrow\Omega_{x}^{0,1}$ in the orthonormal basis $\left\{
\omega_{i}\right\}  $ and $\left\{  \overline{\omega_{i}}\right\}  $ are given
by symmetric matrices. From here $\left(  \ref{LO1}\right)  $ follows
directly. Indeed from the relations of the elements of the matrices of the
operators $\phi_{i}$ in an orthonormal basis $\phi_{i,\overline{\beta}%
}^{\alpha}=\phi_{i,\overline{\alpha}}^{\beta}$ we obtain%
\begin{equation}
\left(  \phi_{i}\circ\overline{\phi_{j}}\right)  _{\overline{\alpha}%
}^{\overline{\beta}}=%
{\displaystyle\sum\limits_{\mu=1}^{n}}
\phi_{i,\overline{\alpha}}^{\mu}\overline{\phi_{j,\overline{\mu}}^{\beta}}=%
{\displaystyle\sum\limits_{\mu=1}^{n}}
\phi_{i,\overline{\mu}}^{\alpha}\overline{\phi_{j,\overline{\beta}}^{\mu}}=%
{\displaystyle\sum\limits_{\mu=1}^{n}}
\overline{\phi_{j,\overline{\beta}}^{\mu}}\phi_{i,\overline{\mu}}^{\alpha
}=\left(  \overline{\phi_{j}}\circ\phi_{i}\right)  _{\beta}^{\alpha}.
\label{LO1a}%
\end{equation}
From the definition of the operator $\left(  \phi_{i}\circ\overline{\phi_{j}%
}\wedge id_{n-2}\right)  ^{\ast}$ given by $\left(  \ref{LO2}\right)  $ we
get:%
\[
\left(  \phi_{i}\circ\overline{\phi_{j}}\wedge id_{n-2}\right)  ^{\ast}\left(
\ast\overline{\omega_{k}}\right)  =\ast\left(  \left(  \phi_{i}\circ
\overline{\phi_{j}}\right)  \left(  \overline{\omega_{k}}\right)  \right)
=\ast\left(
{\displaystyle\sum\limits_{k,l=1}^{n}}
\left(  \phi_{i}\circ\overline{\phi_{j}}\right)  _{k}^{l}\left(
\overline{\omega_{l}}\right)  \right)  =
\]%
\begin{equation}%
{\displaystyle\sum\limits_{k,l=1}^{n}}
\left(  \phi_{i}\circ\overline{\phi_{j}}\right)  _{k}^{l}\left(  \ast\left(
\overline{\omega_{l}}\right)  \right)  =%
{\displaystyle\sum\limits_{k,l=1}^{n}}
\left(  \phi_{i}\circ\overline{\phi_{j}}\right)  _{k}^{l}\left(
\overline{\omega_{1}}\wedge...\wedge\overline{\omega_{l-1}}\wedge
\overline{\omega_{l+1}}\wedge...\wedge\overline{\omega_{n}}\right)
\label{LO3}%
\end{equation}
Combining $\left(  \ref{LO1a}\right)  $ and $\left(  \ref{LO3}\right)  $ we
get%
\[
\left(  \phi_{i}\circ\overline{\phi_{j}}\wedge id_{n-2}\right)  ^{\ast}\left(
\ast\overline{\omega_{k}}\right)  =%
{\displaystyle\sum\limits_{l=1}^{n}}
\left(  \overline{\phi_{j}}\circ\phi_{i}\right)  _{l}^{k}\left(
\overline{\omega_{1}}\wedge...\wedge\overline{\omega_{l-1}}\wedge
\overline{\omega_{l+1}}\wedge...\wedge\overline{\omega_{n}}\right)  =
\]%
\[%
{\displaystyle\sum\limits_{l=1}^{n}}
\left(  \phi_{i}\circ\overline{\phi_{j}}\right)  _{k}^{l}\left(
\overline{\omega_{1}}\wedge...\wedge\overline{\omega_{k-1}}\wedge
\overline{\omega_{k+1}}\wedge...\wedge\overline{\omega_{n}}\right)  .
\]
By using the Calabi-Yau metric and the holomorphic volume form we can identify
the dual of $\Omega_{\text{M}}^{0,n-1}$ with $\Omega_{\text{M}}^{1,0}.$ This
identification is given by
\[
\overline{\omega_{1}}\wedge...\wedge\overline{\omega_{l-1}}\wedge
\overline{\omega_{l+1}}\wedge...\wedge\overline{\omega_{n}}\longrightarrow
\omega_{l}.
\]
Thus from $\left(  \ref{L0}\right)  $ we get
\[
\ast\left(  \left(  \phi_{i}\circ\overline{\phi_{j}}\right)  \left(
\overline{\omega_{k}}\right)  \right)  =
\]%
\begin{equation}
\left(  \phi_{i}\circ\overline{\phi_{j}}\wedge id_{n-2}\right)  ^{\ast}\left(
\ast\overline{\omega_{k}}\right)  =%
{\displaystyle\sum\limits_{l=1}^{n}}
\left(  \overline{\phi_{j}}\circ\phi_{i}\right)  _{l}^{k}\left(  \ast
\overline{\omega_{k}}\right)  =%
{\displaystyle\sum\limits_{l=1}^{n}}
\left(  \overline{\phi_{j}}\circ\phi_{i}\right)  _{l}^{k}\omega_{k}
\label{LO4}%
\end{equation}
From $\left(  \ref{LO4}\right)  $ and $\left(  \ref{LO1a}\right)  $\ we
conclude Proposition \ref{IZS2}. $\blacksquare$

\begin{corollary}
\label{IZS2A}Formula $\left(  \ref{LO4}\right)  $ implies that the composition
of the complex conjugation with the Hodge star operator $\ast$ identifies the
restriction of the image $\left(  \overline{\phi_{j}}\circ\phi_{i}\right)
\left(  \operatorname{Im}\partial\right)  $ on $\operatorname{Im}\partial$ in
$C^{\infty}($M,$\Omega_{\text{M}}^{1,0})$ with the restriction of the image
$\left(  \phi_{i}\circ\overline{\phi_{j}}\wedge id_{n-2}\right)  \left(
\operatorname{Im}\overline{\partial}^{\ast}\right)  $ on $\operatorname{Im}%
\overline{\partial}^{\ast}$ in $C^{\infty}($M,$\Omega_{\text{M}}^{0,n-1})$ and
the restriction of $\overline{\phi_{j}}\circ\phi_{i}\left(  \operatorname{Im}%
\partial^{\ast}\right)  $ on $\operatorname{Im}\partial^{\ast}$ in $C^{\infty
}($M,$\Omega_{\text{M}}^{1,0})$ with the restriction of $\left(  \phi_{i}%
\circ\overline{\phi_{j}}\wedge id_{n-2}\right)  \left(  \operatorname{Im}%
\overline{\partial}\right)  $ on $\operatorname{Im}\overline{\partial}$ in
$C^{\infty}($M,$\Omega_{\text{M}}^{0,n-1}).$
\end{corollary}

From $\left(  \ref{LO4}\right)  $, $\left(  \ref{du}\right)  ,$ and the
identification $\Omega_{\text{M}}^{1,0}$ with $\Omega_{\text{M}}^{0,n-1}$ we
deduce that we can identify $\left(  \overline{\phi_{j}}\circ\phi_{i}\right)
\operatorname{Im}\partial^{\ast}$ in $C^{\infty}($M,$\Omega_{\text{M}}^{1,0})$
with $\left(  \phi_{i}\circ\overline{\phi_{j}}\wedge id_{n-2}\right)
\operatorname{Im}\overline{\partial}$ in $C^{\infty}($M,$\Omega_{\text{M}%
}^{0,n-1})$ and $\left(  \overline{\phi_{j}}\circ\phi_{i}\right)
\operatorname{Im}\partial$ in $C^{\infty}($M,$\Omega_{\text{M}}^{1,0})$ with
$\left(  \phi_{i}\circ\overline{\phi_{j}}\wedge id_{n-2}\right)
\operatorname{Im}\overline{\partial}^{\ast}$ in $C^{\infty}($M,$\Omega
_{\text{M}}^{0,n-1}).$ Since on a K\"{a}hler manifold we have that
$\Delta_{\partial}=\Delta_{\overline{\partial}},$ $\left(  \ref{iz3}\right)  $
is established$.$ Lemma \ref{IZS} is proved. $\blacksquare$

From $\left(  \ref{iz3}\right)  $ we deduce that%
\[
Tr\left(  \exp(-t\Delta_{\tau,1}^{^{\prime}})\mathcal{F}^{^{^{\prime}}}\left(
1,\phi_{i}\circ\overline{\phi_{j}}\right)  +\exp(-t\Delta_{\tau,1}^{^{"}%
})\mathcal{F}^{^{"}}\left(  1,\phi_{i}\circ\overline{\phi_{j}}\right)
\right)  =
\]
\begin{equation}%
{\displaystyle\int\limits_{\text{M}}}
Tr\left(  \exp(-t\Delta_{\tau,1})\left(  \phi_{i}\circ\overline{\phi_{j}%
}\right)  \right)  vol.\label{WP2}%
\end{equation}
We know that
\begin{equation}
Tr\left(  \mathcal{F}^{^{^{\prime}}}\left(  1,\phi_{i}\circ\overline{\phi_{j}%
}\right)  +\mathcal{F}^{^{"}}\left(  1,\phi_{i}\circ\overline{\phi_{j}%
}\right)  \right)  =%
{\displaystyle\int\limits_{\text{M}}}
Tr\left(  \phi_{i}\circ\overline{\phi_{j}}\right)  vol=W.P.\label{WP2a}%
\end{equation}
Combining $\left(  \ref{WP2}\right)  ,$ $\left(  \ref{WP2a}\right)  $ with
Theorem \ref{Main} we deduce that%
\begin{equation}
dd^{c}\left(  \log\det\Delta_{\tau,1}\right)  =dd^{c}\log\left(  \Delta
_{\tau,1}^{^{\prime}}\times\Delta_{\tau,1}^{^{"}}\right)  =-\operatorname{Im}%
W.P.\label{WP3}%
\end{equation}
Theorem \ref{IZS1} is proved. $\blacksquare$

\begin{theorem}
\label{Tr}The relative dualizing sheaf $\pi_{\ast}\left(  \omega
_{\mathcal{X}/\mathfrak{M}_{L}\text{(M)}}\right)  :=\mathcal{L}$ is a trivial
C$^{\infty}$ line bundle.
\end{theorem}

\textbf{Proof: }It is well known fact that a complex line bundle on a complex
manifold $\mathfrak{M}$(M) is topologically trivial if and only if its first
Chern class is zero. According to \cite{To89} the first Chern class of the
relative dualizing sheaf $\pi_{\ast}\left(  \omega_{\mathcal{X}/\mathfrak{M}%
_{L}\text{(M)}}\right)  :=\mathcal{L}$ is the minus the imaginary part of the
Weil-Petersson metric on $\overline{\mathfrak{M}}($M). We observe that the
regularized determinant of the Laplacian of the CY metric is a well defined
function on $\mathfrak{M}($M). Theorem \ref{IZS1} implies that $dd^{c}\left(
\log\det\Delta_{\tau,1}\right)  $ is minus the imaginary part of the
Weil-Petersson metric on $\overline{\mathfrak{M}}($M). Thus the first Chern
class of the relative dualizing sheaf $\pi_{\ast}\left(  \omega_{\mathcal{X}%
/\mathfrak{M}_{L}\text{(M)}}\right)  :=\mathcal{L}$ is represented by the zero
class of cohomology. This proves that the relative dualizing sheaf $\pi_{\ast
}\left(  \omega_{\mathcal{X}/\mathfrak{M}_{L}\text{(M)}}\right)
:=\mathcal{L}$ is topologically trivial on $\mathfrak{M}$(M). Theorem \ref{Tr}
is proved. $\blacksquare$

\section{The Computation on K\"{a}hler Ricci flat Manifold with a Non Trivial
Canonical Bundle}

\begin{theorem}
\label{enr}Suppose that M is $n$ complex dimensional K\"{a}hler Ricci flat
Manifold such that $\left(  K_{\text{M}}\right)  ^{\otimes n}\approxeq
\mathcal{O}_{\text{M}}.$ Then
\begin{equation}
dd^{c}\left(  \log\det\Delta_{\tau,n}\right)  =-\operatorname{Im}W.P.
\label{wp2}%
\end{equation}

\end{theorem}

\textbf{Proof:} The proof of Theorem \ref{enr} is based on Theorem \ref{Main}
which shows that we need to compute
\[
\int\limits_{\text{M}}Tr\left(  \mathcal{F}^{^{\prime}}\left(  n,\left(
\phi_{i}\circ\overline{\phi_{j}}\right)  \right)  )\right)  vol(g).
\]
Here we are using the fact that $\det\Delta_{\tau,n}=\det\Delta_{\tau,n-1}%
^{"}.$ Since we assumed that $H^{n}($M,$\Omega_{\text{M}}^{n}$)$=0$ we get
that
\[
\int\limits_{\text{M}}Tr\left(  \mathcal{F}^{^{\prime}}\left(  n,\left(
\phi_{i}\circ\overline{\phi_{j}}\right)  \right)  )\right)  vol(g)=\int
\limits_{\text{M}}Tr\left(  \left(  \left(  \phi_{i}\circ\overline{\phi_{j}%
}\right)  \wedge id_{n-1}\right)  )\right)  vol(g).
\]
Direct computation using elementary linear algebra shows that if
$A:V\rightarrow V$ is a linear operator acting on an $n$ dimensional vector
space then the trace of the operator $A\wedge id_{n-1}$ acting on $\wedge
^{n}V$ is given by the formula:%
\begin{equation}
Tr(A\wedge id_{n-1})=Tr(A).\label{LA}%
\end{equation}
From formula $\left(  \ref{LA}\right)  $ and the definition of the
Weil-Petersson metric given by $\left(  \ref{wp1}\right)  $ we derive $\left(
\ref{wp2}\right)  .$ Theorem \ref{enr} is proved. $\blacksquare$

\end{document}